\documentclass[11pt,letterpaper,preprint]{imsart}
\usepackage[margin=1in]{geometry}
\usepackage{amsmath,amssymb,amsthm,mathtools,mathrsfs}
\usepackage[numbers]{natbib}
\usepackage{iftex}
\usepackage{todonotes}
\usepackage{comment}
\IfFileExists{bbm.sty}{%
  \usepackage{bbm}
  \newcommand{\Ind}{\mathbbm{1}}
}{%
  \ifPDFTeX
    \newcommand{\Ind}{\mathbf{1}}
  \else
    \usepackage{fontspec}
    \newfontface\IndicatorFont{Latin Modern Math}
    \newcommand{\Ind}{\mathord{\text{\IndicatorFont\symbol{"1D7D9}}}}
  \fi
}
\usepackage[colorlinks=true,allcolors=blue,pdftitle={Sharp normal approximation over convex sets}]{hyperref}
\startlocaldefs

\newcommand{\R}{\mathbf{R}}
\newcommand{\E}{\mathbf{E}}
\renewcommand{\P}{\mathbf{P}}
\newcommand{\tr}{\operatorname{tr}}
\newcommand{\Cl}{\mathscr{C}}
\newcommand{\dc}{\Delta_{\Cl}}
\newcommand{\tv}{d_{\rm TV}}
\newcommand{\inner}[2]{\langle #1,#2\rangle}
\DeclareMathOperator{\Var}{Var}
\DeclareMathOperator{\dist}{dist}
\theoremstyle{plain}
\newtheorem{theorem}{Theorem}[section]
\newtheorem{proposition}[theorem]{Proposition}
\newtheorem{lemma}[theorem]{Lemma}
\newtheoremstyle{boldremark}
  {\medskipamount}{\medskipamount}{\normalfont}{}{\bfseries}{.}{.5em}{}
\theoremstyle{boldremark}
\newtheorem{remark}{Remark}[section]
\numberwithin{equation}{section}
\endlocaldefs
\begin{document}
\pagestyle{plain}
\setlength{\emergencystretch}{1em}
\begin{frontmatter}
\title{Sharp and explicit Berry--Esseen inequalities\\ over convex sets\thanksref{A1}}
\runtitle{Normal approximation over convex sets}
\thankstext{A1}{First appeared on arXiv on Nov. 9, 2026.}

\begin{aug}
\author{\fnms{Mayya} \snm{Zhilova}\thanksref{a2}}
\runauthor{M. Zhilova}

\thankstext[\ensuremath{\ast\ast}]{a2}{Support by the National Science Foundation Award CAREER DMS-2048028 is gratefully acknowledged.} 
\address{School of Mathematics\\ Georgia Institute of Technology\\
Atlanta, GA 30332-0160 USA\\
 e-mail: \href{mailto:mzhilova@math.gatech.edu}{\texttt{mzhilova@math.gatech.edu}}}
\end{aug}

\begin{abstract}
We establish finite-sample Berry--Esseen bounds over all convex subsets of $\R^{d}$
for sums of independent random vectors. The upper bounds give explicit
dependence on the sample size, dimension, and moments of the summands,
with numerical constants. They allow arbitrary dependence among
coordinates and cover every real moment order from three onward.
For i.i.d.\
summands, we determine the optimal dependence of the worst-case
error on the sample size, the dimension, and the prescribed moment
bound. We obtain matching lower bounds using bounded distributions
that can be chosen independently of the sample size.

The results reveal a fourth-moment threshold for normal approximation
over convex sets. Four is the first moment order for which uniform
bounds on the corresponding coordinate moments yield a worst-case
error bound linear in the dimension. Bounds on higher moments give
the same worst-case rate.

The proof combines convexity with the addition of independent Gaussian
noise to establish a normal approximation bound for Rademacher sums
on Parseval frames. Random swaps with independent copies then allow
us to apply this bound conditionally to general sums of independent
random vectors.
\end{abstract}
\begin{keyword}[class=MSC]
\kwd[Primary ]{60F05}
\kwd{62E17}
\kwd[; secondary ]{60E15}
\end{keyword}
\begin{keyword}
\kwd{Multivariate central limit theorem}
\kwd{Berry--Esseen inequality}
\kwd{convex sets}
\kwd{finite sample approximation} 
\kwd{dependence on dimension}
\kwd{anti-concentration}
\kwd{fourth moments}
\end{keyword}
\end{frontmatter}

\section{Introduction and main results}
Estimating the error in the multivariate central limit theorem uniformly
over convex sets is a classical problem in probability. Early work on
this problem includes results of Ranga Rao~\cite{ranga-rao-convex},
Bhattacharya~\cite{bhattacharya-early,bhattacharya-refinements}, Sazonov~\cite{sazonov-convex},
and Nagaev~\cite{nagaev-convex}. Bentkus~\cite{bentkus-dimension} obtained a
third-moment bound with a factor of $d^{1/4}$ and asked whether that
factor could be removed. When the dimension grows with the sample size,
the dependence on dimension and on the moments of the summands
determines whether an error bound guarantees convergence.

Convex-set normal approximation continues to arise in new settings.
Nourdin et al.~\cite{npy-convex} establish bounds over convex sets for
vectors of Wiener--It\^o integrals, and
Schulte and Yukich~\cite{schulte-yukich-geometric} obtain such bounds
for counts of connected components in random geometric graphs.

Gaussian and bootstrap approximations over convex sets help quantify
uncertainty in statistics and machine learning, including in complex
high-dimensional models. Approximations over rectangles and sparsely
convex sets support simultaneous inference in high dimensions, as in
Chernozhukov et al.~\cite{cck-maxima,cck-high-dimensional}. Approximations over Euclidean balls support inference and classical confidence regions for a mean, as in Zhilova~\cite{zhilova-edgeworth,zhilova-bootstrap}. Related approximations enter the analysis of likelihood-based bootstrap confidence sets, including under model misspecification: see \cite{spokoiny-zhilova}, \cite{zhilova-bootstrap}. In machine learning,
Shao and Zhang~\cite{shao-zhang-sgd} quantify the Gaussian approximation
error for parameter estimates from stochastic gradient descent,
and Favaro et al.~\cite{favaro-networks} do so for outputs of
randomly initialized neural networks.

Our results over all convex sets removes the extra
factor $d^{1/4}$ in earlier bounds for sums of independent random
vectors. For four finite moments, Theorem~\ref{t} gives an explicit universal constant below
$16.88$ and allows arbitrary dependence among the coordinates of each
summand. Theorem~\ref{thm:finite-p} gives explicit dimension-independent
bounds for finite moments of every real order $q\in[3,4)$.

For i.i.d.\ sums, Theorem~\ref{thm:moment-constraint} determines the
optimal dependence on the sample size, dimension, and prescribed bound
on the $q$th moment of the Euclidean norm, up to universal constants,
for every real $q\ge3$. The matching lower bounds use simplex and
lattice constructions. For each dimension and moment constraint, a
single bounded law gives the lower bound for all sample sizes. Under
uniform bounds on coordinate fourth moments, the sharp worst-case
error is of order $\min\{1,d/\sqrt n\}$. The simplex example attains
this rate with bounded coordinates, so higher moment assumptions
still give the same worst-case rate. With independent coordinates,
the sharp fourth-moment rate improves to $\min\{1,\sqrt{d/n}\}$.

The upper bounds in our results are explicit in the dimension,
sample size, and moments of the summands, with numerical values or
formulas for the constants.

For the class $\Cl$ of convex Borel subsets of $\R^d$, the convex-set distance is
\[
 \dc(X,Y)=\sup\nolimits_{A\in\Cl}|\P(X\in A)-\P(Y\in A)|\leq 1.
\]
The vector $Z_\Sigma\sim\mathcal N(0,\Sigma)$ is centered Gaussian
with covariance matrix $\Sigma\succeq0$, and $I_d$ denotes the
$d\times d$ identity matrix. A centered vector is called
isotropic if its covariance is $I_d$. The norm $\|\cdot\|$ is Euclidean
for vectors and is the operator norm for matrices. The covariance
matrix of $X$ is denoted by $\Var(X)$.

For summands $X_1,\ldots,X_n$, the normalized sum is
$S_n=n^{-1/2}\sum_{i=1}^n X_i$. For real $q\ge3$, the moment quantities
for an array with total covariance $I_d$ are
$B_q=\sum_{i=1}^n\E\|X_i\|^q$ and $\delta^2=B_4$.
For i.i.d.\ summands, $M_q=\E\|X_1\|^q$. If the i.i.d.\ summands are
centered and isotropic, the rescaled summands $X_i/\sqrt n$ have total
covariance $I_d$.

\subsection{Explicit Berry--Esseen inequalities over convex sets}
\begin{theorem}[Fourth-moment bounds over convex sets]\label{t}
For independent centered random vectors $X_1,\ldots,X_n$ in $\R^d$
with $\sum_{i=1}^n\Var(X_i)=I_d$ and $\delta^2=\sum_{i=1}^n\E\|X_i\|^4<\infty$, the bound is
\begin{equation}\label{eq:main}
 \dc\!\left({\textstyle \sum_{i=1}^nX_i},Z_{I_d}\right)
 \le C_{\mathrm{ind}}\delta,\qquad C_{\mathrm{ind}}<16.88.
\end{equation}
For i.i.d.\ centered isotropic $X_1,\ldots,X_n$, the bound is
\[ \qquad
 \dc(S_n,Z_{I_d})\le \sqrt{M_4/n}\,
 \begin{cases}
 C_{\mathrm{iid},1},&d=1,\\
 C_{\mathrm{iid}},&d\ge2,
 \end{cases}
 \qquad C_{\mathrm{iid},1}<0.95,\quad C_{\mathrm{iid}}<14.81.
\]
Coefficients $C_{\mathrm{ind}},  C_{\mathrm{iid}},  C_{\mathrm{iid,1}}$ are defined in \eqref{eq:main-constants}. For $d\ge2$,
the coefficient $C_{\mathrm{iid}}$ decreases strictly from its value at
$d=2$ to approximately $12.36$
as $d\to\infty$.\\\\
For deterministic $v_1,\ldots,v_n\in\R^d$ satisfying
$\sum_{i=1}^nv_iv_i^T=I_d$ and i.i.d.\ Rademacher random variables
$\varepsilon_1,\ldots,\varepsilon_n$,
\[
 \dc\!\left({\textstyle \sum_{i=1}^n\varepsilon_i v_i},Z_{I_d}\right)
 \le \overline C_{\scriptscriptstyle\mathrm{RS}}\rho,\qquad
 \rho^2=\sum\nolimits_{i=1}^n\|v_i\|^4,\qquad
 \overline C_{\scriptscriptstyle\mathrm{RS}}<10.9 .
\]
The constant $\overline C_{\scriptscriptstyle\mathrm{RS}}$ is defined in
\eqref{eq:rs-uniform-constant}.
\end{theorem}

A family with $\sum_{i=1}^nv_iv_i^T=I_d$ is a \emph{Parseval frame}.
The vector $W=\sum_{i=1}^n\varepsilon_i v_i$ is its Rademacher sum,
where each $\varepsilon_i$ takes $\pm1$ with probability $1/2$.

No assumption is imposed on dependence among the coordinates.
In the i.i.d.\ case, if $\E|X_{1j}|^4\le K_4$ for $1\le j\le d$,
Cauchy--Schwarz and Jensen inequalites give
\[
 d^{2} \leq M_4=\sum\nolimits_{j=1}^d\sum\nolimits_{k=1}^d\E(X_{1j}^2X_{1k}^2)\le K_4d^2 .
\]
Thus $\dc(S_n,Z_{I_d})\le C_{\mathrm{iid}}\sqrt{K_4}\,d/\sqrt n$ for
$d\ge2$, with $C_{\mathrm{iid},1}$ in dimension one.
The simplex example in Theorem~\ref{thm:simplex-lower}(i) with $a=1$
attains the order $\min\{1,d/\sqrt n\}$ with uniformly bounded
coordinates. Thus $d=o(\sqrt n)$ suffices
for uniform convergence under a coordinate fourth-moment bound and
cannot in general be weakened. The coefficients are bounded uniformly in
$d$, and $M_4\ge d^2$ by Jensen's inequality.

For nonsingular total covariance $\Sigma=\sum_{i=1}^n\Var(X_i)\succ0$,
the theorem applies to $\Sigma^{-1/2}X_i$. Affine invariance of $\dc$
then gives the corresponding bound for the original sum.

\begin{theorem}[Finite moments of order below four]\label{thm:finite-p}
For independent centered vectors $X_1,\ldots,X_n$ in $\R^d$ with
$\sum_{i=1}^n\Var(X_i)=I_d$, the bound for real $q$ with $3\le q<4$
and $B_q<\infty$ is
\[
 \dc\!\left({\textstyle \sum_{i=1}^nX_i},Z_{I_d}\right)
 \le C_q^{\mathrm{ind}} B_q^{1/(q-2)}.
\]
For i.i.d.\ centered isotropic $X_1,\ldots,X_n$,
\[
 \dc(S_n,Z_{I_d})\le C_q^{\mathrm{iid}}M_q^{1/(q-2)}/\sqrt n.
\]
Coefficients $C_q^{\mathrm{ind}}, C_q^{\mathrm{iid}} $ are defined in \eqref{eq:finite-moment-constant} and decrease strictly
with $q$ for fixed $d$. The pairs (infimum, maximum) are approximately 
$(28.65,49.63)$ for $C_q^{\mathrm{ind}}$, $(4.68,8.11)$ for
$C_q^{\mathrm{iid}}$ at $d=1$, and $(21.67,44.09)$ for
$C_q^{\mathrm{iid}}$ across $d\ge2$. For $d\ge2$, the i.i.d.\ coefficient
also decreases strictly with $d$ at fixed $q$.
\end{theorem}
\subsection{Sharp finite-sample bounds under a moment constraint}
For integers $d,n\ge1$ and real $q\ge3$, let $\widetilde{M}_q$ be a
prescribed upper bound on $M_q$. For a centered isotropic law, Jensen's
inequality gives $M_q\ge d^{q/2}$, so the prescribed bound satisfies
$d^{q/2}\le\widetilde{M}_q<\infty$. The worst-case error is
\[
 D_q(n,d;\widetilde{M}_q)=\sup\bigl\{\dc(S_n,Z_{I_d}):
   \E X_1=0,\ \Var(X_1)=I_d,\ M_q\le\widetilde{M}_q\bigr\},
\]
where the supremum is over the laws of $X_1$ and $S_n$ is the normalized
sum of $n$ i.i.d.\ copies of $X_1$.
In the following theorem, $A\asymp B$ means $cB\le A\le CB$ for every
admissible $(n,d,q,\widetilde{M}_q)$, with constants $0<c\le C<\infty$
independent of these parameters.

\begin{theorem}[Sharp bounds for i.i.d.\ sums under a moment constraint]\label{thm:moment-constraint}
Uniformly over $d,n,q,\widetilde{M}_q$ in this range,
\begin{equation}\label{eq:moment-constraint}
 D_q(n,d;\widetilde{M}_q)\asymp
 \begin{cases}
 \displaystyle\min\left\{1,
   d/\sqrt{n}+{(\widetilde{M}_q-d^{q/2})^{1/(q-2)}}/{\sqrt n}\right\},&3\le q<4,\\[6pt]
 \displaystyle\min\left\{1,
   {d^{(q-4)/(2(q-2))}\widetilde{M}_q^{1/(q-2)}}/{\sqrt n}\right\},&4\le q<\infty.
 \end{cases}
\end{equation}
For each $q,d,\widetilde{M}_q$,
a single bounded law gives the lower bound for every $n\ge1$.
\end{theorem}

For $3\le q<4$, the excess $\widetilde{M}_q-d^{q/2}$ distinguishes distributions
whose norm is nearly constant from those allowed by a fixed normalized
moment bound. In particular, $\widetilde{M}_q-d^{q/2}=O(d^{q-2})$ already gives
the order $\min\{1,d/\sqrt n\}$.

\subsection{Matching lower bounds}
We give matching lower bounds for i.i.d.\ sums. 
Part (i) below matches the fourth-moment rate. Part (ii) matches the rate for \(3\le q<4\) when the prescribed bound on \(M_q\) is at least \(2d^{q/2}\). 
The lattice and simplex constructions also give the lower bounds
in Theorem~\ref{thm:moment-constraint}.
The proofs are in Appendix~\ref{sec:lower-bounds}.

\begin{theorem}[Lower bounds from a simplex and a lattice example]\label{thm:simplex-lower}
In each part, $X_1,\ldots,X_n$ are independent copies of $X$, and $S_n$
is their normalized sum.

\smallskip
\noindent\textup{(i)}
For every $d\ge1$ and $a\ge1$, there is a bounded centered isotropic law
in $\R^d$ such that, for $X$ with this law and every $q\ge3$,
\[
 \E\|X\|^q=d^{q/2}a^{q-2},\qquad
 \dc(S_n,Z_{I_d})\ge c\min\{1,ad/\sqrt n\},\qquad n\ge1,
\]
where $c>0$ is universal. If $a=1$, the law can be chosen with every
coordinate bounded in absolute value by $\sqrt2$.

\smallskip
\noindent\textup{(ii)}
For every $d\ge1$, $3\le q<4$, and $2d^{q/2}\le\widetilde{M}_q<\infty$,
there is a bounded centered isotropic $X$ in $\R^d$ satisfying
\[
 \E\|X\|^q\le\widetilde{M}_q,\qquad
 \max_{1\le j\le d}\E|X_j|^q\le{\widetilde{M}_q}/{d^{q/2}} .
\]
For this law, simultaneously for all $n\ge1$,
\[
 \dc(S_n,Z_{I_d})\ge\frac3{40}
 \min\left\{1,{\widetilde{M}_q^{1/(q-2)}}/{\sqrt n}\right\}.
\]
In particular, $\widetilde{M}_q=2d^{q/2}$ gives coordinate absolute
$q$th moments at most $2$ and a lower bound of order
$\min\{1,d^{q/[2(q-2)]}/\sqrt n\}$.
\end{theorem}

\begin{remark}\label{rem:f-k-lower}
The lower bound in~\cite{fang-koike-convex} is of order $\sqrt{M_4/n}$
for sufficiently large $n$ and $d\le x_n=\sqrt n/\log n$. The construction uses an
$n$-dependent law and a halfspace perpendicular to the first coordinate axis. 
The atom at $x_n$ contributes $x_n^2/2$ through the fourth moment of that coordinate. The other terms contribute $O(1+d^2)$. Thus $M_4\asymp n/(\log n)^2$ in their
range, and the lower bound has order $1/\log n$.

The example in Theorem~\ref{thm:simplex-lower}(i) with $a=1$
has $M_4=d^2$ and uniformly bounded coordinates. For each dimension,
this fixed law attains the order $\min\{1,\sqrt{M_4/n}\}$ for every
sample size.
\end{remark}

\subsection{Bounds for sums with independent coordinates}
Coordinate independence improves the dimension dependence in the fourth-moment
bound. The upper bound allows the scalar laws to vary across summands and
coordinates. The matching example uses a single scalar law. 

\begin{theorem}[Independent coordinates]\label{thm:independent-coordinates}
Let $X_1,\ldots,X_n$ be independent centered vectors in $\R^d$, each
with independent coordinates and finite fourth moments, and suppose
$n^{-1}\sum_{i=1}^n\Var(X_i)=I_d$. Then
\[
 \dc(S_n,Z_{I_d})\le \frac{C_{\mathrm{ind}}}{n}
 \left(
{\textstyle \sum_{i=1}^n\sum_{j=1}^d\E|X_{ij}|^4}\right)^{1/2},
 \qquad C_{\mathrm{ind}}<16.88.
\]
In particular, if $\E|X_{ij}|^4\le m_4$ for all $1\le i\le n$ and
$1\le j\le d$, the bound is $C_{\mathrm{ind}}\sqrt{dm_4/n}$.

Furthermore, there is a scalar law with
bounded support and a density, mean zero, variance one,
and fourth moment $12/5$ such that, when all $nd$ coordinates have this law
and are independent,
\[
 \dc(S_n,Z_{I_d})\ge c\min\{1,\sqrt{d/n}\}
 \qquad\text{for every }n,d\ge1,
\]
where $c>0$ is universal.
\end{theorem}

The proof is in Appendix~\ref{sec:coordinate-proof}.

\subsection{Discussion of the fourth-moment threshold}\label{sec:fourth-moment-threshold}
For $K>1$ and $\widetilde{M}_q=Kd^{q/2}$, Theorem~\ref{thm:moment-constraint}
gives
\begin{equation}\label{eq:moment-threshold}
 D_q(n,d;Kd^{q/2})\asymp_{K,q}
 \begin{cases}
 \min\{1,d^{q/[2(q-2)]}/\sqrt n\},&3\le q<4,\\
 \min\{1,d/\sqrt n\},&q\ge4,
 \end{cases}
\end{equation}
Here the comparison constants in $\asymp_{K,q}$ may depend only on
$K$ and $q$. The dimension exponent
decreases from $3/2$ at $q=3$ to $1$ at $q=4$ and remains $1$ for
$q>4$.

For $3\le q<4$, the lattice construction in
Section~\ref{sec:lattice-projection} places total mass $L^{-2}$ at
$\pm Lu_1$, where $u_1$ is a unit vector. The two atoms contribute
one unit of variance in direction $u_1$ and $L^{q-2}$ to $M_q$. The bound $M_q\le Kd^{q/2}$ permits
$L\asymp_{K,q} d^{q/[2(q-2)]}$. The projection of $S_n$ then lies on a
lattice of spacing $L/\sqrt n$. A hyperplane containing one lattice
level has zero Gaussian mass and gives the corresponding lower bound.

For $q\ge4$, interpolation gives
$\sqrt{M_4}\le K^{1/(q-2)}d$, so Theorem~\ref{t} gives the upper
bound in \eqref{eq:moment-threshold}. The simplex law in
Theorem~\ref{thm:simplex-lower}(i), with $a=1$, satisfies
$\|X_1\|=\sqrt d$ and $|X_{1j}|\le\sqrt2$ for every coordinate,
yet its convex-set error is at least $c\min\{1,d/\sqrt n\}$.
It attains $M_q=d^{q/2}$ for every $q\ge3$. The lower-bound
proof combines third-order contributions from the simplex directions
using a test with convex sublevel sets. Higher moments can improve
dependence on $K$, but not the worst-case dimension exponent.
At $K=1$, equality in Jensen's inequality forces
$\|X_1\|=\sqrt d$ almost surely. The rate is then
$\min\{1,d/\sqrt n\}$ for every $q\ge3$.

The rates in \eqref{eq:moment-threshold} are also sharp under
dimension-independent coordinate moment bounds. If
$\max_{1\le j\le d}\E|X_{1j}|^q\le K_0$, the triangle inequality in
$L^{q/2}$ gives $M_q\le K_0d^{q/2}$, so the upper bounds follow from
Theorem~\ref{thm:moment-constraint}. For $3\le q<4$, the lattice example
in Theorem~\ref{thm:simplex-lower}(ii) with $\widetilde M_q=2d^{q/2}$
gives the matching lower bound with coordinate $q$th moments at most
$2$. For $q\ge4$, the simplex example above gives the matching lower
bound with uniformly bounded coordinates.

A comparison between third- and fourth-moment bounds also arises in
Wasserstein distance. For probability laws $\mu,\nu$ on $\R^d$ with
finite $s$th moments, the Wasserstein distance of order $s$ is
\[
 W_s(\mu,\nu)=\inf_{\mathcal L(U)=\mu,\ \mathcal L(V)=\nu}
       (\E\|U-V\|^s)^{1/s},\qquad s\ge1.
\]
Here $\mathcal L(U)$ denotes the law of $U$, and the infimum is over
couplings with the indicated marginals. For random vectors,
$W_s(X,Y)=W_s(\mathcal L(X),\mathcal L(Y))$.
Fang et al.~\cite{fsww}
prove $W_2(S_n,Z_{I_d})\le4\sqrt{M_4/n}$. Since $W_1\le W_2$, this also
gives a $W_1$ bound of order $d/\sqrt n$ when $M_4=O(d^2)$.
Fang et al.~\cite{fsww} contrast this with the $d^{3/2}$ dimension
factor in Rai\v{c}'s third-moment $W_1$ bound~\cite{raic-lipschitz}
when $M_3=O(d^{3/2})$, apart from a logarithmic factor.

For $X_1$ uniform on
$\{\pm\sqrt d\,e_j:1\le j\le d\}$, with $e_j$ the $j$th standard basis vector,
Zhai~\cite{zhai-wasserstein} proves
\[
 \liminf_{n\to\infty}\sqrt n\,W_2(S_n,Z_{I_d})\ge d/4,\qquad d\ge1.
\]
This law has $M_q=d^{q/2}$ for every $q\ge3$, in particular $M_4=d^2$.
Thus the factor $d$ in the fourth-moment $W_2$ bound has the correct order.

\subsection{Related work on convex sets and comparison with our results}
For sums of independent centered summands with total covariance $I_d$,
we use the moment quantities $B_3$ and $\delta$ defined above.
For i.i.d.\ sums in fixed dimension, Ranga Rao's fourth-moment
bound over convex sets contains a logarithmic factor in the sample
size~\cite{ranga-rao-convex}. Bhattacharya~\cite{bhattacharya-early}
obtained the $n^{-1/2}$ rate under a finite moment of order $3+\eta$
for some $\eta>0$. Sazonov~\cite{sazonov-convex} and
Bhattacharya~\cite{bhattacharya-refinements} obtained the same
rate under finite third moments. In the latter setting,
Nagaev~\cite{nagaev-convex} and
Bentkus~\cite{bentkus-early} established bounds of orders $dB_3$ and
$d^{1/2}B_3$, respectively. G\"otze~\cite{gotze-convex} used Stein's
method and induction to establish a bound of order $dB_3$ for general
independent summands with finite third moments. Bhattacharya and
Holmes~\cite{bhattacharya-holmes} revisited G\"otze's proof for i.i.d.\
sums with finite third moments and obtained a bound of order $d^{5/2}B_3$.

Bentkus~\cite{bentkus-dimension} subsequently established the bound
$400d^{1/4}B_3$ for i.i.d.\ sums  with finite third moments and extended the order $d^{1/4}B_3$ to
general independent sums~\cite{bentkus-lyapunov}. Rai\v{c}~\cite{raic-explicit}  combined G\"otze's Stein method with Bentkus's smoothing argument and
used spherical symmetry to sharpen the derivative estimates in Bhattacharya
and Holmes's exposition, obtaining the bound $(42d^{1/4}+16)B_3$. In the same
third-moment setting, Ling's recent preprint~\cite{ling-absolute} states the
bound $980B_3$. Our Theorem~\ref{thm:finite-p} yields the bound
$C_3^{\mathrm{ind}}B_3$ with $C_3^{\mathrm{ind}}<49.7$ for three finite moments. For i.i.d.\ sums,
Theorem~\ref{thm:moment-constraint} also identifies
the effect of prescribing a moment close to its minimum. When
$\widetilde{M}_3-d^{3/2}=O(d)$, the sharp worst-case error is of order
$\min\{1,d/\sqrt n\}$.

Zhai~\cite{zhai-wasserstein} bounds $\dc(T,Z_{I_d})$ by
$5d^{1/6}W_2(T,Z_{I_d})^{2/3}$. For independent centered isotropic
summands with $\|X_i\|\le\beta$ almost surely, this comparison and
his $W_2$ estimate give a convex-set bound of order
$d^{1/2}\beta^{2/3}(1+\log n)^{2/3}n^{-1/3}$.
Under the same assumptions, Theorem~\ref{t} gives
$16.88\beta\sqrt{d/n}$, since
$\E\|X_i\|^4\le\beta^2d$.

For independent summands with finite fourth moments and total
covariance $I_d$, the bound in~\cite{fang-koike-convex} is
$Cd^{1/4}\bar\delta|\log\bar\delta|$, where
$\bar\delta=\min\{\delta,e^{-1}\}$. Cai et al.~\cite{cai-yang-zhang} study statistics
that decompose into a sum of independent vectors and a remainder. For independent
summands with finite fourth moments, their bound $Cd^{1/4}\delta$
removes the logarithmic factor from the bound in~\cite{fang-koike-convex}.
Under the same moment and covariance assumptions,
Theorem~\ref{t} removes the additional factor $d^{1/4}$ and obtains
$C_{\mathrm{ind}}\delta$ with $C_{\mathrm{ind}}<16.88$.

Remark~\ref{rem:f-k-lower} compares the lower bound in~\cite{fang-koike-convex}
with our lower bound, which holds for a fixed bounded law in each
dimension and every sample size. Theorem~\ref{thm:moment-constraint}
determines, up to universal constants, the worst-case error for
i.i.d.\ sums under each prescribed moment bound.

Kasprzak and Peccati~\cite{kasprzak-peccati} establish convex-distance
bounds for vector-valued functions of independent inputs in terms of
moments of resampling differences. For isotropic i.i.d.\ sums with
finite sixth moments, their estimate has the $n^{-1/2}$ rate and an
explicit factor $d^4$.

Nourdin et al.~\cite{npy-convex} establish convex-distance bounds
for centered, Malliavin-differentiable functionals of a Gaussian
field with a positive-definite Gaussian target covariance. Their results include
fourth-moment bounds for vectors of multiple Wiener--It\^o integrals.
Favaro et al.~\cite{favaro-networks} treat vectors with square-integrable
Stein kernels, as well as centered
Gaussian mixtures with square-integrable random covariance. Our
bounds for independent sums require only the stated moment and
covariance assumptions. For the example of Favaro et al.\ with
Gaussian initialization, one input, and one hidden layer,
Theorem~\ref{t} yields $O(d/\sqrt n)$ under a fixed bounded activation
and nondegenerate output covariance, where $n$ is the hidden width
and $d$ is the output dimension. Favaro et al. obtain $O(d^{7/4}/\sqrt n)$ for this example.

\subsection{Proof outline and organization of the paper}\label{sec:map}

\paragraph*{Fourth-moment bounds.}
Section~\ref{sec:rademacher} establishes a dimension-free fourth-moment
bound over convex sets for Rademacher sums on deterministic Parseval
frames (Proposition~\ref{prop:frame}). The two one-sided bounds use
Gaussian coupling of two blocks and smoothing with covariance adapted
to the frame, respectively (Propositions~\ref{prop:minus}
and~\ref{prop:plus}). In Section~\ref{sec:independent}, a conditional
application of the frame bound after random swaps with independent
copies compares the original sum with a Gaussian-smoothed average
of the sum and its independent copy (Proposition~\ref{prop:swap}).
The independent Gaussian component permits a Taylor expansion that
controls the remaining error in total variation
(Proposition~\ref{prop:fixed}). Section~\ref{sec:main-proof} combines
the two comparisons to prove Theorem~\ref{t}.

\paragraph*{Other moment orders.}
Section~\ref{sec:proof-finite-p} proves Theorem~\ref{thm:finite-p}
using a cutoff on the summand norms. A conditional application of
Theorem~\ref{t} gives a Gaussian comparison with fixed covariance.
Convolution with this Gaussian permits a Taylor comparison of the
remaining terms. For the upper bounds in Theorem~\ref{thm:moment-constraint},
Section~\ref{sec:proof-moment-upper} bounds the tail $q$th moment in
terms of $M_q-d^{q/2}$ when $3\le q<4$ and applies
Lemma~\ref{lem:moment-cutoff}. For $q\ge4$, Theorem~\ref{t} applies with the
fourth moment bounded by H\"older's inequality when $q>4$.

\paragraph*{Lower bounds.}
In Appendix~\ref{sec:lattice-projection}, an atom in a one-dimensional
lattice projection gives the lower bound in
Theorem~\ref{thm:simplex-lower}(ii).
For the simplex example in Theorem~\ref{thm:simplex-lower}(i),
Appendix~\ref{sec:proof-simplex} combines a lattice argument with
a Taylor expansion for a test function whose sublevel sets are
convex. Appendix~\ref{sec:proof-moment-lower} chooses between the
lattice and simplex constructions to obtain the lower bounds in
Theorem~\ref{thm:moment-constraint}.

\paragraph*{Independent coordinates.}
In Appendix~\ref{sec:coordinate-proof}, the upper bound in
Theorem~\ref{thm:independent-coordinates} follows from
Theorem~\ref{t} applied to the individual coordinate contributions.
For the lower bound, a bounded scalar law with nonzero third moment
gives a cubic contribution in a coordinatewise Gaussian replacement.

Appendix~\ref{sec:gaussian-perimeter-proof} proves the Gaussian
perimeter bound in Lemma~\ref{lem:shell}, which controls the error from
removing the smoothing. Appendix~\ref{sec:scalar-tilt-proof} proves
the scalar integral estimate in Lemma~\ref{lem:scalar-tilt}, which
enters the density comparison in Lemma~\ref{lem:density}.

\section*{Notation}
The notation $e_j$ denotes the $j$th standard basis vector and
$\mathbf 1$ the vector of ones.
For a class \(\mathscr A\) of Borel sets, \(\Delta_{\mathscr A}\) is the supremum over \(\mathscr A\) of the absolute difference of probabilities. Total variation is the case of all Borel sets. The distance notation
applies to both laws and random vectors.
For $\Sigma\succ0$, $\varphi_\Sigma:\R^d\to(0,\infty)$ is the density
of $Z_\Sigma$. The functions $\varphi_1$ and $\Phi$ are the standard
univariate normal density and distribution function, respectively.

For symmetric matrices, $A\succeq B$ means that $A-B$ is positive
semidefinite, and $A\succ B$ means that it is positive definite.
For $A\succeq0$, $A^{1/2}$ is the positive semidefinite square root.
The transpose of $A$ is denoted by $A^T$, its trace by
$\tr A$, and the Euclidean adjoint of a linear map $B$ by $B^*$.
For arrays of the same shape, $\inner{\cdot}{\cdot}$ is the entrywise
inner product and $\|\cdot\|_{\mathrm F}$ is the corresponding Frobenius norm.
For $h\in\R^d$, $h^{\otimes k}$ is the $k$-fold tensor product of $h$.
If $f:\R^d\to\R$ is $k$ times differentiable, $f^{(k)}(x)$ is its $k$th
derivative tensor, with
$(h^T\nabla)^kf(x)=\inner{f^{(k)}(x)}{h^{\otimes k}}$, $x,h\in\R^d$.

For a Borel set $A$, $\Ind_A:\R^d\to\{0,1\}$ is its indicator.
The symbol $*$ denotes convolution, and $\|g\|_\infty$ denotes the
supremum of $|g|$ on its domain. The positive part is
$x_+=\max\{x,0\}$, $x\in\R$.

%
\section{Gaussian approximation for Rademacher sums over convex sets}\label{sec:rademacher}
Throughout this section, $W=\sum_{i=1}^n\varepsilon_i v_i$, where
$\varepsilon_1,\ldots,\varepsilon_n$ are i.i.d.\ Rademacher random variables,
$\sum_{i=1}^n v_iv_i^T=I_d$, and
$\rho^2\coloneqq\sum_{i=1}^n\|v_i\|^4$. For each convex set $A$, we bound
$\P(Z_{I_d}\in A)-\P(W\in A)$ and $\P(W\in A)-\P(Z_{I_d}\in A)$ separately.

\subsection{Gaussian coupling and the first one-sided bound}
The first bound couples two block sums to independent Gaussians using a
quadratic Wasserstein estimate. Convexity bounds
$\P(Z_{I_d}\in A)-\P(W\in A)$ by translation errors, each controlled
using the Gaussian from the other block.

The distance $W_2$ is the case $s=2$ of the Wasserstein distance
defined in Section~\ref{sec:fourth-moment-threshold}. For an optimal
coupling $(U,V)$ of $\mu,\nu$ in $W_2$, Cauchy--Schwarz gives
$\E\|U-V\|\le W_2(\mu,\nu)$. This expected displacement enters the
translation comparison.

\begin{lemma}[Quadratic Wasserstein bound for a Rademacher sum]\label{lem:transport}
For $v_1,\ldots,v_m\in\R^d$ satisfying $\sum_{i=1}^{m} v_iv_i^T=I_d$ and
i.i.d.\ Rademacher random variables $\varepsilon_1,\ldots,\varepsilon_m$,
the sum $R=\sum_{i=1}^{m}\varepsilon_i v_i$ and the quantity
$r^2=\sum_{i=1}^{m}\|v_i\|^4$ satisfy
\[
 W_2(R,Z_{I_d})\le C_{W_2}r,\qquad C_{W_2}=4.
\]
\end{lemma}
\begin{proof}
The estimate of Fang et al.~\cite{fsww} implies that independent centered
vectors $T_1,\ldots,T_m$ with finite fourth moments and total covariance
$I_d$, and an independent copy $(\tilde T_i)_{i=1}^m$, satisfy
\begin{equation}\label{eq:fsww}
 W_2^2\left({\textstyle \sum_{i=1}^{m} T_i},Z_{I_d}\right)\le 2\sum\nolimits_{i=1}^{m}\E\|T_i-\tilde T_i\|^4 .
\end{equation}
Applying their estimate to $X_i=\sqrt m\,T_i$, squaring, and using the
second term in each minimum bounds the left-hand side of
\eqref{eq:fsww} by
\[
 \frac{2}{m^2}\sum\nolimits_{i=1}^{m}\sum\nolimits_{j=1}^d
 \E\big[\|X_i-\tilde X_i\|^2(X_{ij}-\tilde X_{ij})^2\big]
 =2\sum\nolimits_{i=1}^{m}\E\|T_i-\tilde T_i\|^4.
\]
The equality uses $X_i-\tilde X_i=\sqrt m\,(T_i-\tilde T_i)$ and
sums the squared coordinate differences over $j$.

For $T_i=\varepsilon_i v_i$, the frame identity gives total covariance
$I_d$. For an independent copy $\tilde\varepsilon_i$,
$|\varepsilon_i-\tilde\varepsilon_i|$ equals $0$ or $2$, each with
probability $1/2$. Thus $\E\|T_i-\tilde T_i\|^4=8\|v_i\|^4$, and
\eqref{eq:fsww} gives $W_2^2(R,Z_{I_d})\le C_{W_2}^2r^2$ with $C_{W_2}^2=16$.
\end{proof}

\begin{lemma}[Partition into two blocks and Gaussian couplings]\label{lem:partition}
For $0<\rho<1$, the constants are $a=(1-\rho)/2$ and $a_+=(1+\rho)/2$.
There is a partition $\{1,\ldots,n\}=s_1\sqcup s_2$ such that, with
\[
 R_{s_j}={\textstyle\sum_{i\in s_j}\varepsilon_i v_i},\qquad
 \Sigma_{s_j}=\Var(R_{s_j}),\qquad
 \tau_{s_j}^2={\textstyle\sum_{i\in s_j}\|v_i\|^4},\qquad j\in\{1,2\},
\]
one has
\[
 aI_d\preceq\Sigma_{s_j}\preceq a_+I_d,\qquad \tau_{s_1}\le \tau_{s_2},\qquad
 \tau_{s_1}^2+\tau_{s_2}^2=\rho^2.
\]
For each normalized block,
\[
 r_{s_j}=\left({\textstyle\sum_{i\in s_j}\|\Sigma_{s_j}^{-1/2}v_i\|^4}\right)^{1/2}
 \le \tau_{s_j}/a,\qquad j\in\{1,2\}.
\]
There are independent couplings
$(R_{s_j},G_{s_j})$ with $G_{s_j}\sim\mathcal{N}(0,\Sigma_{s_j})$ and
\begin{equation}\label{cost}
 \E\|R_{s_j}-G_{s_j}\|\le C_{W_2}\sqrt{a_+}\,\tau_{s_j}/a.
\end{equation}
\end{lemma}
\begin{proof}
Each index in $\{1,\ldots,n\}$ is assigned independently to the first or
second block with probability $1/2$. With $\xi_i$ the indicator of
assignment to the first block, the two block covariance matrices are
\[
 \Sigma_1=\sum\nolimits_{i=1}^{n} \xi_i v_iv_i^T,\qquad \Sigma_2=I_d-\Sigma_1.
\]
Since $\Sigma_1-I_d/2=\sum_{i=1}^{n}(\xi_i-1/2)v_iv_i^T$, independence and
centering make all cross terms vanish in the expected squared
Frobenius norm. Thus
\[
 \E\|\Sigma_1-I_d/2\|_{\mathrm F}^2
 =\frac14\sum\nolimits_{i=1}^{n}\|v_iv_i^T\|_{\mathrm F}^2=\frac{\rho^2}{4}.
\]
Some assignment therefore has squared deviation at most $\rho^2/4$.
For such an assignment,
the operator norm of $\Sigma_1-I_d/2$ is at most its Frobenius norm, and
$\Sigma_2-I_d/2=-(\Sigma_1-I_d/2)$. Both covariance matrices therefore
have eigenvalues in $[a,a_+]$.

Since $a>0$, both blocks are nonempty. They are labeled $s_1,s_2$ so that
$\tau_{s_1}\le \tau_{s_2}$. The blocks partition $\{1,\ldots,n\}$, so
$\tau_{s_1}^2+\tau_{s_2}^2=\rho^2$.

For each $j\in\{1,2\}$, the normalized vectors
$w_i=\Sigma_{s_j}^{-1/2}v_i$, $i\in s_j$, satisfy
\[
 {\textstyle\sum_{i\in s_j}w_iw_i^T=\Sigma_{s_j}^{-1/2}\Sigma_{s_j}\Sigma_{s_j}^{-1/2}=I_d},
 \qquad
 r_{s_j}=\left({\textstyle \sum_{i\in s_j}\|w_i\|^4}\right)^{1/2}\le \tau_{s_j}/a.
\]
The inequality follows from $\|\Sigma_{s_j}^{-1/2}\|\le a^{-1/2}$.
By Lemma~\ref{lem:transport},
\[
 W_2\!\left(\mathcal L(\Sigma_{s_j}^{-1/2}R_{s_j}),\mathcal N(0,I_d)\right)
 \le C_{W_2}r_{s_j} .
\]
An optimal quadratic coupling $(U_{s_j},Z_{s_j})$ of these laws satisfies
\[
 \E\|U_{s_j}-Z_{s_j}\|^2
 = W_2^2\!\left(\mathcal L(\Sigma_{s_j}^{-1/2}R_{s_j}),\mathcal N(0,I_d)\right)
 \le C_{W_2}^2r_{s_j}^2.
\]
With $\widehat R_{s_j}=\Sigma_{s_j}^{1/2}U_{s_j}$ and $G_{s_j}=\Sigma_{s_j}^{1/2}Z_{s_j}$,
the vector $\widehat R_{s_j}$ has the law of $R_{s_j}$ and
$G_{s_j}\sim\mathcal N(0,\Sigma_{s_j})$. By Cauchy--Schwarz and
$\|\Sigma_{s_j}^{1/2}\|\le\sqrt{a_+}$,
\[
 \E\|\widehat R_{s_j}-G_{s_j}\|
 \le \bigl(\E\|\widehat R_{s_j}-G_{s_j}\|^2\bigr)^{1/2}
 \le C_{W_2}\sqrt{a_+}\,r_{s_j}
 \le C_{W_2}\sqrt{a_+}\,\tau_{s_j}/a .
\]

The resulting pairs $(\widehat R_{s_1},G_{s_1})$ and
$(\widehat R_{s_2},G_{s_2})$ are sampled independently. For
$j\in\{1,2\}$, the possible sign assignments for a given sum form the set
\[
 \mathcal S_{s_j}(x)=
 \left\{e\in\{-1,1\}^{s_j}:{\textstyle \sum_{i\in s_j}e_i v_i}=x\right\},
 \qquad x\in\R^d.
\]
In block $s_j$, an auxiliary random variable uniform on $[0,1]$ selects the sign vector
$\boldsymbol\varepsilon_{s_j}=(\varepsilon_i)_{i\in s_j}$ uniformly from
$\mathcal S_{s_j}(\widehat R_{s_j})$. The auxiliary variables are independent
of each other and of both coupled pairs. The set is nonempty almost surely
because $\widehat R_{s_j}$ has the law of the block sum.

For any sign assignment $e$ with $x=\sum_{i\in s_j}e_i v_i$, the
construction can select $e$ only when $\widehat R_{s_j}=x$.
For i.i.d.\ Rademacher random variables, exactly $|\mathcal S_{s_j}(x)|$ of the
$2^{|s_j|}$ equally likely assignments give the sum $x$. Thus
\[
 \begin{aligned}
 \P(\boldsymbol\varepsilon_{s_j}=e)
 &=\P(\widehat R_{s_j}=x)\,
   \P(\boldsymbol\varepsilon_{s_j}=e\mid \widehat R_{s_j}=x)\\
 &=\frac{2^{-|s_j|}|\mathcal S_{s_j}(x)|}{|\mathcal S_{s_j}(x)|}
 =2^{-|s_j|}.
 \end{aligned}
\]
This proves that the Rademacher random variables within each block are i.i.d.\
The two block constructions use independent inputs, so
$(\boldsymbol\varepsilon_{s_1},\widehat R_{s_1},G_{s_1})$ and
$(\boldsymbol\varepsilon_{s_2},\widehat R_{s_2},G_{s_2})$ are independent.
In particular, each $G_{s_j}$ is independent of the entire sign vector
in the other block.

The block sum $R_{s_j}$ formed from the selected signs equals
$\widehat R_{s_j}$ almost surely. Hence the pairs $(R_{s_j},G_{s_j})$ are
independent and satisfy \eqref{cost}.
\end{proof}

\begin{lemma}[Comparison under translation]\label{lem:slice}
The vector $X$ has density $p_X$ on $\R^d$, and $\Sigma$ is a symmetric
$d\times d$ matrix with $\Sigma\succeq aI_d$ for some $a>0$. The density
assumption is that, for some $M\ge1$,
\[
 p_X(x)\le M\varphi_\Sigma(x),\qquad\text{for almost every }x\in\R^d,
\]
where $\varphi_\Sigma$ is the density of the centered Gaussian with covariance $\Sigma$.
For every convex set $A$ and $h\in\R^d$,
\[
 \max\{\P(X\in A,\ X+h\notin A),\ \P(X\notin A,\ X+h\in A)\}
 \le M\varphi_1(0)\|h\|/\sqrt a .
\]
The same bound holds for $|\P(X+h\in A)-\P(X\in A)|$.
\end{lemma}
\begin{proof}
If $h=0$, both events are empty and the difference in probabilities is
zero. The remaining case is $h\ne0$. The entering and leaving sets are
$C_+=\{x:x\notin A,\ x+h\in A\}$ and
$C_-=\{x:x\in A,\ x+h\notin A\}$, respectively.

With $t=\|h\|$ and $e=h/t$, every point has a unique representation
$z+ue$ with $z\in e^\perp$ and $u\in\R$. By convexity,
$I_z=\{u\in\R:z+ue\in A\}$ is an interval, possibly empty or unbounded.
Translation by $h$ changes $u$ to $u+t$ on the same line.

If $u\in I_z$ and $u+t\notin I_z$, no $v\in I_z$ can satisfy
$v>u+t$. Otherwise $u+t$ would lie between $u$ and $v$ and hence
belong to $I_z$. Thus $b_z=\sup I_z$ is finite and
$u\le b_z\le u+t$. The leaving section is therefore empty or
contained in $[b_z-t,b_z]$.

If $u\notin I_z$ and $u+t\in I_z$, no $v\in I_z$ can satisfy
$v<u$, since then $u$ would lie between $v$ and $u+t$ and hence
belong to $I_z$. Thus $a_z=\inf I_z$ is finite and
$u\le a_z\le u+t$. The entering section is therefore empty or
contained in $[a_z-t,a_z]$.

Each section of either crossing set has length at most $t=\|h\|$,
regardless of whether $I_z$ contains its endpoints.

For $C\in\{C_+,C_-\}$, its section is
$C_z=\{u:z+ue\in C\}$. In an orthonormal basis
$e_1,\ldots,e_{d-1}$ of $e^\perp$, the vector $z$ has the representation
$z=\sum_{j=1}^{d-1}z_j e_j$. The map
$(z_1,\ldots,z_{d-1},u)\mapsto z+ue$ has columns
$e_1,\ldots,e_{d-1},e$ and absolute determinant one. Hence
$dx=dz\,du$, where $dz$ is Lebesgue measure on $e^\perp$ and $du$
is length measure. For $d=1$, $dz$ is unit mass at $0$.

For $G\sim\mathcal N(0,\Sigma)$, the decomposition is $G=Z+Ue$, with
$U=e^TG$ and $Z=G-Ue$. The joint density of $(Z,U)$ is
$\varphi_\Sigma(z+ue)$, so the marginal density of $Z$ is
$m(z)=\int_\R\varphi_\Sigma(z+ve)\,dv$, $z\in e^\perp$.
With $Q=\Sigma^{-1}$, $c=e^TQe>0$, and $b(z)=e^TQz$ for $z\in e^\perp$,
completing the square gives
\[
 (z+ue)^TQ(z+ue)
 =z^TQz+2b(z)u+cu^2
 =c\bigl(u+b(z)/c\bigr)^2+z^TQz-b(z)^2/c.
\]
Thus $\varphi_\Sigma(z+ue)=K(z)\exp[-c(u+b(z)/c)^2/2]$ for a
positive function $K:e^\perp\to(0,\infty)$ independent of $u$. Substituting
$w=\sqrt c\,(v+b(z)/c)$ in the integral defining $m(z)$ gives
$m(z)=K(z)\sqrt{2\pi/c}$. The conditional density of $U$ given
$Z=z$ is therefore
\[
 \frac{\varphi_\Sigma(z+ue)}{m(z)}
 =\sqrt{\frac{c}{2\pi}}\,
   \exp\!\left[-\frac c2\bigl(u+b(z)/c\bigr)^2\right],
 \qquad z\in e^\perp,\ u\in\R.
\]
Since $\Sigma\succeq aI_d$, we have $Q\preceq a^{-1}I_d$ and
$c\le a^{-1}$. The conditional density is bounded by
$\varphi_1(0)/\sqrt a$, uniformly in $z$. By Fubini and
$\int_{e^\perp}m(z)\,dz=1$,
\[
 \P(G\in C)
 =\int_{e^\perp}m(z)\int_{C_z}
   \frac{\varphi_\Sigma(z+ue)}{m(z)}\,du\,dz
 \le \frac{t\varphi_1(0)}{\sqrt a}.
\]
Integrating $p_X\le M\varphi_\Sigma$ over $C$ now gives
$\P(X\in C)\le M t\varphi_1(0)/\sqrt a$.

The event that both $X$ and $X+h$ belong to $A$ contributes equally
to the two probabilities.
Cancelling that term gives
\[
 \P(X+h\in A)-\P(X\in A)=\P(X\in C_+)-\P(X\in C_-).
\]
The bound for both crossing probabilities therefore gives
\[
 \bigl|\P(X+h\in A)-\P(X\in A)\bigr|
 \le\max\{\P(X\in C_+),\P(X\in C_-)\}
 \le {M\varphi_1(0)\|h\|}/{\sqrt a}.
\]
\end{proof}

\begin{proposition}[One-sided bound by Gaussian coupling]\label{prop:minus}
For $0<\rho<1$, with $a=(1-\rho)/2$ and $a_+=(1+\rho)/2$, every convex
set $A$ satisfies
\[
 \P(Z_{I_d}\in A)-\P(W\in A)
 \le 2\sqrt2\,C_{W_2}\varphi_1(0)\frac{\sqrt{a_+}}{a^{3/2}}\,\rho
 =\frac{4C_{W_2}}{\sqrt\pi}\frac{\sqrt{1+\rho}}{(1-\rho)^{3/2}}\,\rho.
\]
\end{proposition}
\begin{proof}
By Lemma~\ref{lem:partition}, there are independent coupled pairs
$(R_{s_1},G_{s_1})$ and $(R_{s_2},G_{s_2})$, with
$R_{s_1}+R_{s_2}=W$ and $G_{s_j}\sim\mathcal N(0,\Sigma_{s_j})$.
The sum and block differences are $G=G_{s_1}+G_{s_2}$ and
$D_{s_j}=R_{s_j}-G_{s_j}$.
The Gaussian vectors are independent and their covariances add to
$I_d$, so $G\sim\mathcal N(0,I_d)$. Also,
\[
 W=G+D_{s_1}+D_{s_2}
   =\tfrac12(G+2D_{s_1})+\tfrac12(G+2D_{s_2}).
\]
If both $G+2D_{s_1}$ and $G+2D_{s_2}$ belong to $A$, then their
midpoint $W$ belongs to $A$ by convexity. Thus
\[
 \{G\in A,\ W\notin A\}\subseteq
 \{G\in A,\ G+2D_{s_1}\notin A\}
 \cup
 \{G\in A,\ G+2D_{s_2}\notin A\}.
\]
Cancelling the common term $\P(G\in A,\ W\in A)$ and applying the
union bound gives
\[
 \begin{aligned}
 \P(Z_{I_d}\in A)-\P(W\in A)
 &=\P(G\in A,\ W\notin A)-\P(G\notin A,\ W\in A)\\
 &\le{\textstyle\sum_{j=1}^2\P(G\in A,\ G+2D_{s_j}\notin A)}.
 \end{aligned}
\]
For $j\in\{1,2\}$ and $k=3-j$, the crossing event is
$E_j=\{G\in A,\ G+2D_{s_j}\notin A\}$. Conditional on
$(R_{s_j},G_{s_j})=(r,g)$, the two points are $G_{s_k}+g$ and
$G_{s_k}+2r-g$. They share the independent Gaussian
$G_{s_k}\sim\mathcal N(0,\Sigma_{s_k})$ and differ by $2(r-g)$.
The function $F_j:\R^d\times\R^d\to[0,1]$ giving the conditional crossing probability is
\[
 F_j(r,g)=\P\bigl(G_{s_k}\in A-g,\
                  G_{s_k}+2(r-g)\notin A-g\bigr),
 \qquad r,g\in\R^d.
\]
The set $A-g$ is convex, and $\Sigma_{s_k}\succeq aI_d$.
Lemma~\ref{lem:slice}, with $X=G_{s_k}$, $M=1$, and $h=2(r-g)$, gives
\[
 F_j(r,g)\le\frac{2\varphi_1(0)}{\sqrt a}\|r-g\|,\qquad r,g\in\R^d.
\]
By independence, evaluating $F_j$ at the random pair gives a version
of the conditional probability of $E_j$ given $(R_{s_j},G_{s_j})$.
Consequently,
\[
 \P(E_j\mid R_{s_j},G_{s_j})
 =F_j(R_{s_j},G_{s_j})
 \le\frac{2\varphi_1(0)}{\sqrt a}\|D_{s_j}\|
 \qquad\text{almost surely}.
\]
Taking expectations yields
\[
 \P(E_j)\le\frac{2\varphi_1(0)}{\sqrt a}\E\|D_{s_j}\|.
\]

The union bound and the estimate for $\P(E_j)$ give
\[
 \P(Z_{I_d}\in A)-\P(W\in A)
 \le\sum\nolimits_{j=1}^2\P(E_j)
 \le\frac{2\varphi_1(0)}{\sqrt a}
       \sum\nolimits_{j=1}^2\E\|D_{s_j}\|.
\]
By \eqref{cost},
$\E\|D_{s_j}\|\le C_{W_2}\sqrt{a_+}\,\tau_{s_j}/a$. Therefore
\[
 \P(Z_{I_d}\in A)-\P(W\in A)
 \le 2C_{W_2}\varphi_1(0)\frac{\sqrt{a_+}}{a^{3/2}}
       (\tau_{s_1}+\tau_{s_2}).
\]
The blocks partition the indices, so $\tau_{s_1}^2+\tau_{s_2}^2=\rho^2$.
Cauchy--Schwarz yields $\tau_{s_1}+\tau_{s_2}\le\sqrt2\,\rho$, giving the
stated bound.
\end{proof}

\subsection{Density estimates and smoothing with an independent Gaussian component}
To bound $\P(W\in A)-\P(Z_{I_d}\in A)$, we first add independent Gaussian
noise. A scalar integral estimate gives density bounds for the smoothed
sum and the smoothed smaller block. The density bound for the smaller
block controls the replacement of the other block by a Gaussian. The
bound for the full smoothed sum is used to remove the noise.

Here $\log_+ x=\max\{0,\log x\}$ for $x>0$, and the smoothing parameters are
\begin{gather}
 r_*=0.22,\qquad A_0=1.1,\qquad L=0.26, \label{eq:smoothing-parameters}\\
 \lambda_r(\beta)=\beta\bigl(A_0+L\log_+(r^2/\beta)\bigr),
 \qquad 0<\beta\le r\le r_*. \label{eq:smoothing-lambda}
\end{gather}

For a Parseval frame $v_1,\ldots,v_m$ of nonzero vectors, the smoothing
parameters use $\beta_i=\|v_i\|^2$, $1\le i\le m$, and
$r=(\sum_{i=1}^m\beta_i^2)^{1/2}$.
The matrix $B=\sum_{i=1}^m\lambda_r(\beta_i)v_iv_i^T$ is the covariance of the
centered Gaussian vector $Z_B\sim N(0,B)$ that we add independently to the
associated Rademacher sum. This sum has covariance $I_d$, so adding $Z_B$
gives covariance $I_d+B$. Adding this noise gives the discrete sum a smooth
density, which we compare pointwise with the density of $N(0,I_d+B)$.
The constants $A_0$ and $L$ determine how much noise is added along each
frame direction $v_i$. They are chosen to make the density comparison
possible while keeping $\|B\|_{\mathrm F}$ small enough to control the
error when the noise is removed.

\begin{lemma}[A scalar integral under exponential tilting]\label{lem:scalar-tilt}
For $0<\beta\le r\le r_*$ and $\lambda=\lambda_r(\beta)$, the integral is
\[
 I(s)=\frac1{\sqrt{2\pi}}\int_{\R} e^{-\lambda t^2/2}
       \bigl(1-\operatorname{sech}^2 s\,\sin^2(\sqrt\beta\,t)\bigr)^{1/(2\beta)}\,dt,\qquad s\in\R,
\]
where $\operatorname{sech}s=1/\cosh s$. The function $K$ is
\[
 K(r)=\max\left\{\frac1{16}+\frac{35}{4}r^2,\;12r^2\right\},
 \qquad 0<r\le r_*.
\]
Then $0<I(s)<\infty$ and
\begin{equation}\label{eq:app-a-scalar}
 \log\cosh s+\beta\log\!\bigl(I(s)\sqrt{1+\lambda}\bigr)
 \le \frac{s^2}{2}+\frac{\beta^2}{r^2}K(r),\qquad s\in\R.
\end{equation}
\end{lemma}

The proof of Lemma~\ref{lem:scalar-tilt} is given in Appendix~\ref{sec:scalar-tilt-proof}.

\begin{lemma}[Density bound for a smoothed Rademacher sum]\label{lem:density}
The vectors $v_1,\ldots,v_m\in\R^d$ are nonzero and satisfy
$\sum_{i=1}^{m} v_iv_i^T=I_d$. With $\beta_i=\|v_i\|^2$ and
$r^2=\sum_{i=1}^{m}\beta_i^2$, the parameters for $r\le r_*$ are
\[
 \lambda_i=\lambda_r(\beta_i),\qquad
 B=\sum\nolimits_{i=1}^{m}\lambda_i v_iv_i^T,\qquad \Sigma=I_d+B .
\]
Then $B\succ0$ and
\[
 \|B\|_{\mathrm F}\le (A_0+2L/e)r^2<1.292r^2.
\]
For i.i.d.\ Rademacher random variables $\varepsilon_1,\ldots,\varepsilon_m$
and an independent $Z_B\sim\mathcal N(0,B)$, the density $p$ of
$Y=\sum_{i=1}^{m}\varepsilon_i v_i+Z_B$ satisfies
\begin{equation}\label{eq:app-a-density}
 p(y)\le C_{\mathrm{pdf}}(r)\varphi_\Sigma(y),\qquad y\in\R^d,
\end{equation}
where
\[
 C_{\mathrm{pdf}}(r)
 =\exp\!\left(K(r)+\frac{A_0^2}{4}r^3\right)<2,
 \qquad 0<r\le r_*,
\]
and $K(r)=\max\{1/16+(35/4)r^2,\,12r^2\}$ as in Lemma~\ref{lem:scalar-tilt}.
\end{lemma}

\begin{proof}
The frame identity gives
\[
 B\succeq(\min_{1\le i\le m}\lambda_i)I_d\succ0.
\]
The covariance also has the decomposition
\[
 B=A_0\sum_{i=1}^{m}\beta_i v_iv_i^T+V^T\operatorname{diag}(b_i)V,\qquad
 b_i=L\beta_i\log_+(r^2/\beta_i),
\]
where $V$ has rows $v_i^T$. The first term is positive semidefinite
with trace $A_0r^2$, which bounds its Frobenius norm.
Since $V^TV=I_d$,
\[
 \|V^T\operatorname{diag}(b_i)V\|_{\mathrm F}^2
 \le \|\operatorname{diag}(b_i)V\|_{\mathrm F}^2
 =\sum\nolimits_{i=1}^{m}\beta_i b_i^2
 \le \frac{4L^2}{e^2}r^4 .
\]
The last bound uses $\sqrt\beta\log_+(r^2/\beta)\le 2r/e$, obtained
for $\beta\le r^2$ by maximizing $u e^{-u/2}$ with
$u=\log(r^2/\beta)$. Hence $\|B\|_{\mathrm F}\le(A_0+2L/e)r^2$.
The stated parameters give $A_0+2L/e<1.292$.

To compare $p$ with $\varphi_\Sigma$, we apply Fourier inversion to
$p_\theta(x)=e^{\theta^Tx}p(x)$, $x\in\R^d$, for $\theta\in\R^d$.
Here $s_i=\theta^Tv_i$, $1\le i\le m$.
Since $p$ is a finite Gaussian mixture,
$p_\theta$ is continuous and integrable. With the Fourier kernel
$e^{i\xi^Tx}$, independence gives
\[
 \widehat p_\theta(\xi)
 =\exp\!\left(\tfrac12(\theta+i\xi)^TB(\theta+i\xi)\right)
   \prod\nolimits_{i=1}^{m}\cosh(s_i+i v_i^T\xi),\qquad \xi\in\R^d.
\]
This transform is integrable by $B\succ0$ and
$|\cosh(s+it)|\le\cosh s$. Fourier inversion gives
\[
 p(y)\le(2\pi)^{-d}e^{-\theta^Ty+\theta^TB\theta/2}
 \int_{\R^d} e^{-\xi^TB\xi/2}\prod\nolimits_{i=1}^{m}|\cosh(s_i+i v_i^T\xi)|\,d\xi,
 \qquad y\in\R^d.
\]
The integrand has a product form to which the geometric Brascamp--Lieb
inequality applies, reducing the integral to one-dimensional factors
estimated by Lemma~\ref{lem:scalar-tilt}.
For $u_i=v_i/\sqrt{\beta_i}$ and $p_i=\beta_i$, the linear maps are
$B_i:\R^d\to\R$, $B_i\xi=u_i^T\xi$ for $\xi\in\R^d$, and the functions
$f_i:\R\to[0,\infty)$ are given by
\[
 f_i(t)=e^{-\lambda_i t^2/2}|\cosh(s_i+i\sqrt{\beta_i}\,t)|^{1/\beta_i},
 \qquad t\in\R.
\]
The integrand is $\prod_{i=1}^{m} f_i(B_i\xi)^{p_i}$.
Since $B_i^*z=zu_i$ and $\|u_i\|=1$, the maps and weights satisfy
\[
 B_iB_i^*=I_1,\qquad
{\textstyle \sum_{i=1}^{m} p_i B_i^*B_i=\sum\nolimits_{i=1}^{m}\beta_i u_i u_i^T=I_d}.
\]
The geometric Brascamp--Lieb inequality, in the formulation of
Bennett et al.~\cite{bcct}, bounds the integral by
$\prod_{i=1}^{m}(\int_{\R} f_i)^{p_i}$. Each $f_i$ is nonnegative and integrable. The identity
$|\cosh(s+it)|^2=\cosh^2 s\,(1-\operatorname{sech}^2s\sin^2t)$ gives
$\int_\R f_i=\sqrt{2\pi}\cosh(s_i)^{1/\beta_i}I_i(s_i)$, where $I_i$
is the function on $\R$ from Lemma~\ref{lem:scalar-tilt}
at $(\beta_i,\lambda_i)$.
Since $\sum_{i=1}^{m}\beta_i=d$ by the frame identity,
\begin{equation}\label{eq:app-a-product}
 p(y)\le (2\pi)^{-d/2}e^{-\theta^Ty+\theta^TB\theta/2}
                 \prod\nolimits_{i=1}^{m}\bigl[\cosh(s_i)I_i(s_i)^{\beta_i}\bigr],\qquad y\in\R^d.
\end{equation}

Since $0<\beta_i\le(\sum_{j=1}^{m}\beta_j^2)^{1/2}=r\le r_*$,
Lemma~\ref{lem:scalar-tilt} gives
\[
 \log\cosh s_i+\beta_i\log I_i(s_i)
 \le \frac{s_i^2}{2}+\frac{\beta_i^2}{r^2}K(r)
       -\frac{\beta_i}{2}\log(1+\lambda_i).
\]
Together with $\sum_{i=1}^{m}s_i^2=\|\theta\|^2$ and
$\sum_{i=1}^{m}\beta_i^2=r^2$, this reduces \eqref{eq:app-a-product} to
\[
 p(y)\le (2\pi)^{-d/2}
 \exp\!\left(-\theta^Ty+\frac12\theta^T\Sigma\theta+K(r)
             -\frac12{\textstyle\sum_{i=1}^{m}\beta_i\log(1+\lambda_i)}\right),
 \qquad y\in\R^d.
\]
Taking $\theta=\Sigma^{-1}y$ and comparing with the density
$\varphi_\Sigma$ gives
\[
 p(y)\le e^{K(r)+D_B}\varphi_\Sigma(y),\qquad
 D_B=\frac12\left[\log\det(I_d+B)
              -{\textstyle\sum_{i=1}^{m}\beta_i\log(1+\lambda_i)}\right].
\]
Since $\log\det(I_d+B)\le\tr B=\sum_{i=1}^{m}\beta_i\lambda_i$ and
$x-\log(1+x)=\int_0^x t/(1+t)\,dt\le x^2/2$ for $x\ge0$,
\[
 D_B\le\frac14\sum\nolimits_{i=1}^{m}\beta_i\lambda_i^2.
\]
For $u=\log_+(r^2/\beta_i)$, the inequality $A_0+Lu\le A_0e^{u/2}$
follows from $L\le A_0/2$ and $e^{u/2}\ge1+u/2$. Hence
\[
 \frac{\lambda_i}{\sqrt{\beta_i}}
 \le A_0\sqrt{\beta_i}\,e^{u/2}
 =A_0\max\{\sqrt{\beta_i},r\}\le A_0\sqrt r.
\]
It follows that $D_B\le(A_0^2r/4)\sum_{i=1}^{m}\beta_i^2=A_0^2r^3/4$,
which proves \eqref{eq:app-a-density}.
The exponent of $C_{\mathrm{pdf}}(r)$ increases with $r$ and is at
most $K(r_*)+A_0^2r_*^3/4<\log 2$.
\end{proof}

\begin{lemma}[Derivatives of Gaussian-smoothed indicators]
\label{lem:gaussian-convolution}
Let $A\subseteq\R^d$ be Borel, let $\Sigma\succ0$, and set
$f_A=\Ind_A*\varphi_\Sigma$. For an order-$k$ tensor $U$,
let $\operatorname{Sym}U$ be the average of $U$ over all $k!$
permutations of its indices. For every integer $k\ge1$ and $x\in\R^d$,
\[
 \left|\inner{U}{f_A^{(k)}(x)}\right|
 \le \frac{\sqrt{k!}}{2}
 \left\|(\Sigma^{-1/2})^{\otimes k}\operatorname{Sym}U\right\|_{\mathrm F},
\]
where the tensor power acts by applying $\Sigma^{-1/2}$ to each tensor
factor. In particular, for $v\in\R^d$,
\[
 |(v^T\nabla)^k f_A(x)|
 \le \frac{\sqrt{k!}}{2}(v^T\Sigma^{-1}v)^{k/2},
 \qquad
 \|f_A^{(k)}(x)\|_{\mathrm F}
 \le \frac{\sqrt{k!}}{2}\|\Sigma^{-1/2}\|^k.
\]
\end{lemma}
\begin{proof}
Let $K_U(y)=\inner{U}{\varphi_\Sigma^{(k)}(y)}$ and
$V=(\Sigma^{-1/2})^{\otimes k}\operatorname{Sym}U$.
The symmetry of the derivative tensor and the change of variables
$y=\Sigma^{1/2}z$ give
\[
 \int_{\R^d}|K_U(y)|\,dy
 =\int_{\R^d}\big|\inner{V}{\varphi_{I_d}^{(k)}(z)}\big|\,dz
 \le\sqrt{k!}\,\|V\|_{\mathrm F}.
\]
The last inequality follows from Cauchy--Schwarz and the Hermite
$L^2$ estimate in the proof of Lemma~9.1 of the supplement to
\cite{zhilova-edgeworth}. A directional form of the Gaussian density
estimate appears in Lemma~A.5 of the supplement to
\cite{zhilova-bootstrap}.

Gaussian density derivatives are integrable, and $K_U$ integrates
to zero since $k\ge1$. Differentiating the convolution gives
\[
 \left|\inner{U}{f_A^{(k)}(x)}\right|
 =\left|\int_{\R^d}(\Ind_A(x-y)-1/2)K_U(y)\,dy\right|
 \le\frac12\int_{\R^d}|K_U(y)|\,dy .
\]
Taking $U=v^{\otimes k}$ gives the directional bound. The Frobenius
bound follows by duality, since
$\|V\|_{\mathrm F}\le\|\Sigma^{-1/2}\|^k\|U\|_{\mathrm F}$.
\end{proof}

For the smoothing comparison, we use the cutoff $\rho_*=1/8$.
For $0<\rho\le\rho_*$, the density estimate of Lemma~\ref{lem:density}
applies to both the full frame and the normalized smaller block, since
$\rho_*<r_*$ and $\sqrt{2}\rho_*/(1-\rho_*)=\sqrt{2}/7<r_*$.

\begin{proposition}[Gaussian approximation after smoothing]\label{prop:smooth}
For $0<\rho\le\rho_*$, with $a=(1-\rho)/2$ and $a_+=(1+\rho)/2$,
the bounds use
\[
 \begin{aligned}
 C_B(\rho)&=(A_0+2L/e)\left(1+\frac{a_+}{2a^2}\right),\\
 \eta_\rho
 &= C_{W_2}\varphi_1(0)\frac{\sqrt{a_+}}{a^{3/2}}\,\rho
    \sup_{0<z\le1/2}\left[C_{\mathrm{pdf}}\!\left(\frac{\rho\sqrt z}{a}\right)
           \sqrt{1-z}\right]
    +\frac{\rho^2}{2\sqrt6\,a^2}.
 \end{aligned}
\]
Here $A_0$ and $L$ are the parameters in the definition of $\lambda_r$ above,
and $C_{\mathrm{pdf}}$ is given by Lemma~\ref{lem:density}.
There exists a positive-definite covariance matrix $B$ such that, for
$Z_B\sim\mathcal N(0,B)$ independent of $W$, the vector $Y=W+Z_B$ satisfies
\begin{equation}\label{eta}
 \begin{aligned}
 \|B\|_{\mathrm F}&\le C_B(\rho)\rho^2,\qquad
 \dc(Y,Z_{I_d+B})\le\eta_\rho,\\
 p_Y(y)&\le C_{\mathrm{pdf}}(\rho)\varphi_{I_d+B}(y),\qquad y\in\R^d.
 \end{aligned}
\end{equation}
On this range, $C_B(\rho)<3.211$, $\eta_\rho<5.625\rho$, and
$C_{\mathrm{pdf}}(\rho)<1.374$.
\end{proposition}
\begin{proof}
We construct a single $Z_B$ that gives density bounds for both
$W+Z_B$ and $R_{s_1}+Z_B$. The latter bound allows us to replace the
other block using Lemma~\ref{lem:slice}.

Lemma~\ref{lem:partition} gives a partition
$s_1\sqcup s_2=\{1,\ldots,n\}$ with
\[
 aI_d\preceq\Sigma_{s_j}\preceq a_+I_d,\qquad
 \tau_{s_1}\le \tau_{s_2},\qquad \tau_{s_1}^2+\tau_{s_2}^2=\rho^2,
\]
with $a$ and $a_+$ as in the statement. Hence $\tau_{s_1}^2\le\rho^2/2$.

Zero coefficient vectors remain zero after normalization and may be
omitted in both applications of Lemma~\ref{lem:density}, since they
contribute zero to the Rademacher sum, its covariance, and the sum
of fourth powers.

For the full frame, $\rho\le\rho_*<r_*$. Lemma~\ref{lem:density} gives
a positive-definite smoothing
covariance $B_g$ with
$\|B_g\|_{\mathrm F}\le(A_0+2L/e)\rho^2$.

For the small block, Lemma~\ref{lem:partition} shows that the vectors
$w_i=\Sigma_{s_1}^{-1/2}v_i$, $i\in s_1$, form a Parseval frame and
that $r_{s_1}\le\tau_{s_1}/a$. Hence
\[
 r_{s_1}\le\frac{\rho}{\sqrt2\,a}
 =\frac{\sqrt2\,\rho}{1-\rho}
 \le\frac{\sqrt2\,\rho_*}{1-\rho_*}<r_* .
\]
Thus Lemma~\ref{lem:density} gives a positive-definite smoothing
covariance $\widetilde B_{s_1}$ for the normalized block, with
$\|\widetilde B_{s_1}\|_{\mathrm F}\le(A_0+2L/e)r_{s_1}^2$.

Multiplication by $\Sigma_{s_1}^{1/2}$ returns the normalized block
and its smoothing Gaussian to the original coordinates. The transformed
Gaussian has covariance
\[
 B_{s_1}=\Sigma_{s_1}^{1/2}\widetilde B_{s_1}\Sigma_{s_1}^{1/2}\succ0 .
\]
For compatible matrices, $\|HC\|_{\mathrm F}\le\|H\|\|C\|_{\mathrm F}$
by applying the operator norm bound to each column of $C$.
Applying this inequality on both sides and using
$\|\Sigma_{s_1}^{1/2}\|^2=\|\Sigma_{s_1}\|\le a_+$ gives
\[
 \|B_{s_1}\|_{\mathrm F}
 \le a_+\|\widetilde B_{s_1}\|_{\mathrm F}
 \le(A_0+2L/e)a_+r_{s_1}^2
 \le(A_0+2L/e)a_+\tau_{s_1}^2/a^2.
\]

The Gaussians $Z_{B_g}\sim\mathcal N(0,B_g)$ and
$Z_{B_{s_1}}\sim\mathcal N(0,B_{s_1})$ can be chosen independently
of each other and of both block couplings. Their sum
$Z_B=Z_{B_g}+Z_{B_{s_1}}$ is independent of $W$ and has law
$\mathcal N(0,B)$ with $B=B_g+B_{s_1}\succ0$.
By the triangle inequality,
\[
 \begin{aligned}
 \|B\|_{\mathrm F}
 &\le(A_0+2L/e)\left(\rho^2+\frac{a_+\tau_{s_1}^2}{a^2}\right)
 \le C_B(\rho)\rho^2.
 \end{aligned}
\]
For the density $p_g$ of $W+Z_{B_g}$, Lemma~\ref{lem:density}
gives the density multiplier $C_{\mathrm{pdf}}(\rho)$ for the full frame,
so
\[
 p_g(y)\le C_{\mathrm{pdf}}(\rho)\varphi_{I_d+B_g}(y),\qquad y\in\R^d.
\]
Since $Z_{B_{s_1}}$ is independent of $W+Z_{B_g}$, the density of
$Y=W+Z_{B_g}+Z_{B_{s_1}}$ satisfies
\[
 \begin{aligned}
 p_Y(y)
 &=\int_{\R^d}p_g(y-z)\varphi_{B_{s_1}}(z)\,dz\\
 &\le C_{\mathrm{pdf}}(\rho)\int_{\R^d}
       \varphi_{I_d+B_g}(y-z)\varphi_{B_{s_1}}(z)\,dz\\
 &=C_{\mathrm{pdf}}(\rho)\varphi_{I_d+B}(y),\qquad y\in\R^d.
 \end{aligned}
\]
The last equality uses $B=B_g+B_{s_1}$ and the covariance formula
for a sum of independent centered Gaussians.

For the small block, $H=\Sigma_{s_1}^{1/2}$ and
$\widetilde Z=H^{-1}Z_{B_{s_1}}$. The vector $\widetilde Z$ is
independent of $R_{s_1}$ and has law
$\mathcal N(0,\widetilde B_{s_1})$ by the definition of $B_{s_1}$.
The normalized smoothed sum $\widetilde Y_1=H^{-1}R_{s_1}+\widetilde Z$
has density
\[
 \widetilde p_1(x)=
 \E\!\left[\varphi_{\widetilde B_{s_1}}\bigl(x-H^{-1}R_{s_1}\bigr)\right],
 \qquad x\in\R^d.
\]
Indeed, given $R_{s_1}$, the law of $\widetilde Y_1$ is Gaussian with
mean $H^{-1}R_{s_1}$ and covariance $\widetilde B_{s_1}$.
Lemma~\ref{lem:density}, applied to the normalized block, gives
\[
 \widetilde p_1(x)
 \le C_{\mathrm{pdf}}(r_{s_1})\varphi_{I_d+\widetilde B_{s_1}}(x),
 \qquad x\in\R^d.
\]
Since $H\widetilde Y_1=R_{s_1}+Z_{B_{s_1}}$, its density $p_1$ is
obtained by the change of variables $x\mapsto H^{-1}x$. Thus
\[
 \begin{aligned}
 p_1(x)
 &=(\det H)^{-1}\widetilde p_1(H^{-1}x)\\
 &\le C_{\mathrm{pdf}}(r_{s_1})(\det H)^{-1}
       \varphi_{I_d+\widetilde B_{s_1}}(H^{-1}x)\\
 &=C_{\mathrm{pdf}}(r_{s_1})\varphi_{\Sigma_{s_1}+B_{s_1}}(x),
 \qquad x\in\R^d.
 \end{aligned}
\]
Here $H(I_d+\widetilde B_{s_1})H^T=\Sigma_{s_1}+B_{s_1}$.
Adding the independent Gaussian $Z_{B_g}$ now gives
\[
 \begin{aligned}
 p_{R_{s_1}+Z_B}(y)
 &=\int_{\R^d}p_1(y-z)\varphi_{B_g}(z)\,dz\\
 &\le C_{\mathrm{pdf}}(r_{s_1})\varphi_{\Sigma_{s_1}+B}(y),
 \qquad y\in\R^d.
 \end{aligned}
\]
The density bound for $U=R_{s_1}+Z_B$ allows us to replace the second
block by its Gaussian. The original and modified smoothed sums are
$Y=U+R_{s_2}$ and $Y'=U+G_{s_2}$.
The vector $U$ is independent of the pair $(R_{s_2},G_{s_2})$.
The density of $U$ is bounded by
$C_{\mathrm{pdf}}(r_{s_1})\varphi_{\Sigma_{s_1}+B}$, with
$\Sigma_{s_1}+B\succeq aI_d$.

For a convex set $A$, the shifted probabilities are
\[
 K_A(v)=\P(U+v\in A),\qquad v\in\R^d.
\]
By independence, conditioning on $(R_{s_2},G_{s_2})$ gives
\[
 \begin{aligned}
 \P(Y\in A\mid R_{s_2},G_{s_2})&=K_A(R_{s_2}),\\
 \P(Y'\in A\mid R_{s_2},G_{s_2})&=K_A(G_{s_2})
 \end{aligned}
 \qquad\text{almost surely}.
\]
For fixed $r,g\in\R^d$, their difference evaluated at the pair $(r,g)$ is
\[
 K_A(r)-K_A(g)
 =\P\bigl(U+(r-g)\in A-g\bigr)-\P(U\in A-g),
 \qquad r,g\in\R^d.
\]
The set $A-g$ is convex. Lemma~\ref{lem:slice} with
$X=U$, $M=C_{\mathrm{pdf}}(r_{s_1})$, and $h=r-g$ therefore gives
\[
 |K_A(r)-K_A(g)|
 \le \frac{C_{\mathrm{pdf}}(r_{s_1})\varphi_1(0)}{\sqrt a}\,\|r-g\|,
 \qquad r,g\in\R^d.
\]
The conditional probability difference is
$V=\P(Y\in A\mid R_{s_2},G_{s_2})-\P(Y'\in A\mid R_{s_2},G_{s_2})$.
The law of total expectation, the conditional bound, and \eqref{cost} yield
\[
 \begin{aligned}
 \bigl|\P(Y\in A)-\P(Y'\in A)\bigr|
 &=|\E V|\le\E|V|\\
 &\le\frac{C_{\mathrm{pdf}}(r_{s_1})\varphi_1(0)}{\sqrt a}\,
        \E\|R_{s_2}-G_{s_2}\|\\
 &\le C_{W_2}\varphi_1(0)\frac{\sqrt{a_+}}{a^{3/2}}\,
        C_{\mathrm{pdf}}(r_{s_1})\tau_{s_2}.
 \end{aligned}
\]

Since $\tau_{s_1}\le \tau_{s_2}$ and $\tau_{s_1}^2+\tau_{s_2}^2=\rho^2$,
the ratio $z=\tau_{s_1}^2/\rho^2$ satisfies $0<z\le1/2$ and
$\tau_{s_2}=\rho\sqrt{1-z}$. The strict lower bound follows from
$\Sigma_{s_1}\succeq aI_d$ with $a>0$. The definition in
Lemma~\ref{lem:density} makes $C_{\mathrm{pdf}}(r)$ nondecreasing in $r$.
Since $r_{s_1}\le \tau_{s_1}/a=\rho\sqrt z/a\le r_*$, the last
estimate is uniform over $A$ and gives
\[
 \dc(Y,Y')\le C_{W_2}\varphi_1(0)\frac{\sqrt{a_+}}{a^{3/2}}\,\rho
 \sup_{0<z\le1/2}\left[C_{\mathrm{pdf}}\!\left(\frac{\rho\sqrt z}{a}\right)
             \sqrt{1-z}\right].
\]

The Gaussian $G_{s_2}+Z_B$, independent of the Rademacher random
variables indexed by $s_1$, supplies the smoothing for the Taylor
estimates on the first block. Its covariance is
$\Sigma=\Sigma_{s_2}+B\succeq aI_d$. For a convex set $A$, the smoothed
indicator $f:\R^d\to[0,1]$ is
\[
 f(x)=(\Ind_A*\varphi_\Sigma)(x)
     =\P(x+Z_\Sigma\in A),\qquad x\in\R^d.
\]
Then
\[
 \P(Y'\in A)=\E f(R_{s_1}).
\]
Here $\partial_v^k f(x)=(v^T\nabla)^k f(x)$ for
$x,v\in\R^d$ and $k=0,\ldots,4$. Lemma~\ref{lem:gaussian-convolution}
and $\Sigma\succeq aI_d$ give
\[
 \|\partial_v^4f\|_\infty
 \le\frac{\sqrt{4!}}2(v^T\Sigma^{-1}v)^2
 \le\frac{\sqrt{4!}}{2a^2}\|v\|^4,\qquad v\in\R^d.
\]

For an enumeration $s_1=\{i_1,\ldots,i_m\}$ and i.i.d.\ $\mathcal N(0,1)$
random variables $(\gamma_i)_{i\in s_1}$ independent of all variables
already constructed, the successive replacement sums are
\[
 R^{(k)}=\sum\nolimits_{\ell=1}^{k}\gamma_{i_\ell}v_{i_\ell}
       +\sum\nolimits_{\ell=k+1}^{m}\varepsilon_{i_\ell}v_{i_\ell},
 \qquad k=0,\ldots,m,
\]
where an empty sum is zero. The sum of the unchanged terms in the $k$th
replacement is
\[
 T_k=\sum\nolimits_{\ell=1}^{k-1}\gamma_{i_\ell}v_{i_\ell}
       +\sum\nolimits_{\ell=k+1}^{m}\varepsilon_{i_\ell}v_{i_\ell},
 \qquad k=1,\ldots,m.
\]
Then $R^{(k-1)}=T_k+\varepsilon_{i_k}v_{i_k}$ and
$R^{(k)}=T_k+\gamma_{i_k}v_{i_k}$, with $T_k$ independent of both scalars.
Taylor's formula gives
\[
 f(x+tv)=\sum_{j=0}^{3}\frac{t^j}{j!}\partial_v^j f(x)
          +\mathcal R_{x,v}(t),
 \qquad x,v\in\R^d,\ t\in\R,
\]
where
\[
 |\mathcal R_{x,v}(t)|
 \le\frac{|t|^4}{4!}\|\partial_v^4f\|_\infty,
 \qquad x,v\in\R^d,\ t\in\R.
\]
The scalars $\varepsilon_{i_k}$ and $\gamma_{i_k}$ have matching first
three moments and fourth moments $1$ and $3$, respectively. The Taylor
polynomials therefore have equal expectations. For every $x\in\R^d$,
\[
 \begin{aligned}
 &\bigl|\E f(x+\varepsilon_{i_k}v_{i_k})
           -\E f(x+\gamma_{i_k}v_{i_k})\bigr|\\
 &\qquad\le\frac{1+3}{4!}\|\partial_{v_{i_k}}^4f\|_\infty
 \le\frac{\|v_{i_k}\|^4}{\sqrt6\,a^2}.
 \end{aligned}
\]
Conditioning on $T_k$ and then averaging gives the same bound on
$|\E f(R^{(k-1)})-\E f(R^{(k)})|$. Summing over $k=1,\ldots,m$ gives
\[
 \begin{aligned}
 |\E f(R^{(0)})-\E f(R^{(m)})|
 &\le\sum\nolimits_{k=1}^{m}|\E f(R^{(k-1)})-\E f(R^{(k)})|\\
 &\le\frac{1}{\sqrt6\,a^2}\sum\nolimits_{i\in s_1}\|v_i\|^4
 =\frac{\tau_{s_1}^2}{\sqrt6\,a^2}.
 \end{aligned}
\]
Here $R^{(0)}=R_{s_1}$. The vectors $R^{(m)}$, $G_{s_2}$, and $Z_B$
are independent centered Gaussians with covariances $\Sigma_{s_1}$,
$\Sigma_{s_2}$, and $B$, respectively. Since $s_1$ and $s_2$ partition
$\{1,\ldots,n\}$,
\[
 \Sigma_{s_1}+\Sigma_{s_2}
 =\sum\nolimits_{i=1}^{n}v_iv_i^T=I_d.
\]
Consequently,
\[
 R^{(m)}+G_{s_2}+Z_B\sim\mathcal N(0,I_d+B),\qquad
 \E f(R^{(m)})=\P(Z_{I_d+B}\in A).
\]
Using $\tau_{s_1}^2\le\rho^2/2$ gives
\[
 |\P(Y'\in A)-\P(Z_{I_d+B}\in A)|
 \le\frac{\tau_{s_1}^2}{\sqrt6\,a^2}
 \le\frac{\rho^2}{2\sqrt6\,a^2}.
\]
The triangle inequality and the two replacement bounds give
\[
 \dc(Y,Z_{I_d+B})
 \le \dc(Y,Y')+\frac{\rho^2}{2\sqrt6\,a^2}\le\eta_\rho .
\]
Together with the covariance and density bounds above, this proves
\eqref{eta}.

The factor $A_0+2L/e$ comes from Lemma~\ref{lem:density}.
On $0<\rho\le\rho_*$, the bounds $a\ge7/16$ and $a_+\le9/16$ and
the definition of $C_B(\rho)$ give $C_B(\rho)<3.211$.
For $0<z\le1/2$, the argument $r=\rho\sqrt z/a$ satisfies $r\le r_*$ and
$r^2\le (2\rho_*/(1-\rho_*))^2z=4z/49$. Hence
\[
 \log C_{\mathrm{pdf}}(r)
 \le \frac1{16}+\left(12+\frac{A_0^2r_*}{4}\right)r^2
 \le\frac1{16}+z.
\]
The derivative of $e^z\sqrt{1-z}$ is $e^z(1-2z)/(2\sqrt{1-z})\ge0$
for $0<z<1/2$. Thus the supremum in $\eta_\rho$ is at most
$e^{9/16}/\sqrt2$. Using $\rho\le\rho_*$ in the last term gives
\[
 \frac{\eta_\rho}{\rho}
 \le\frac{24C_{W_2}e^{9/16}}{7\sqrt{7\pi}}+\frac{16}{49\sqrt6}.
\]
For $C_{W_2}=4$, the right-hand side is less than $5.625$.
Monotonicity gives $C_{\mathrm{pdf}}(\rho)\le
C_{\mathrm{pdf}}(\rho_*)<1.374$.
\end{proof}

\subsection{Dimension-free Gaussian perimeter bound}
For a Gaussian measure $\gamma=\mathcal N(0,\Sigma)$ with symmetric
positive definite covariance matrix $\Sigma\in\R^{d\times d}$,
Nazarov~\cite{nazarov} defines the maximal Gaussian perimeter by
\[
 \Gamma(\gamma)
 =\sup
     \left\{\frac{\gamma(C_t\setminus C)}{t}: C\in\Cl, t>0 \right\},
 \qquad
 C_t=\{x\in\R^d:\dist(x,C)\le t\},\qquad t\ge0.
\]
For $x\in\R^d$ and $D\subseteq\R^d$, the distance is
$\dist(x,D)=\inf_{y\in D}\|x-y\|$, with $\dist(x,\varnothing)=+\infty$.
The Euclidean unit ball is $B_2^d=\{x\in\R^d:\|x\|\le1\}$.

For the standard Gaussian, Ball~\cite{ball-gsa} proved the upper bound
$4d^{1/4}$, and Nazarov~\cite{nazarov} showed that the order
$d^{1/4}$ is sharp. Bentkus~\cite{bentkus-dimension,bentkus-lyapunov}
used Gaussian perimeter bounds in convex-set Berry--Esseen estimates.
Rai\v{c}~\cite{raic-explicit} refined the upper bound
to $\sqrt{2/\pi}+0.59(d^{1/4}-1)$. Subsequent work on the lower bound
includes \cite{nadimpalli-pascale,ivanisvili-xie}.

For $\Sigma\succ0$, Nazarov~\cite{nazarov} bounded
$\Gamma(\mathcal N(0,\Sigma))$ above and below by universal positive
constants times $\|\Sigma^{-1}\|_{\mathrm F}^{1/2}$. We use an explicit
upper constant to control the error from removing the smoothing. The proof in
Appendix~\ref{sec:gaussian-perimeter-proof} uses Nazarov's polar
coordinate system, the divergence theorem, and
the variance of a Gaussian quadratic form.

\begin{lemma}[An explicit constant in Nazarov's dimension-free bound for maximal Gaussian perimeter]\label{lem:shell}
For a symmetric positive definite covariance matrix $\Sigma$ and
$\gamma=\mathcal N(0,\Sigma)$, every convex set $C\subseteq\R^d$ satisfies
\begin{equation}\label{shell}
 \gamma(C_t\setminus C)
 \le C_{\scriptscriptstyle\mathrm{GSA}}\,t\|\Sigma^{-1}\|_{\mathrm F}^{1/2},
 \qquad t\ge0.
\end{equation}
Here $C_{\scriptscriptstyle\mathrm{GSA}}=(1+\sqrt3)/\sqrt{2\pi}<1.1$.
\end{lemma}
\subsection{Removing the smoothing and combining the bounds}
Proposition~\ref{prop:smooth} compares $W+Z_B$ with $Z_{I_d+B}$.
The perimeter bound in Lemma~\ref{lem:shell} controls the error of
removing $Z_B$. The resulting bound on $\P(W\in A)-\P(Z_{I_d}\in A)$,
combined with Proposition~\ref{prop:minus}, gives the convex-set distance.

\begin{lemma}[Comparison with the unsmoothed law]\label{lem:remove}
The smoothing vector $Z_B\sim\mathcal{N}(0,B)$ is independent of $W$,
with $B\succ0$ and $b=\|B\|_{\mathrm{F}}\le1/2$. If $Y=W+Z_B$ and
$\Sigma=I_d+B$ satisfy $p_Y(y)\le M\varphi_\Sigma(y)$ for almost every
$y\in\R^d$ for some $M\ge1$ and $\dc(Y,Z_\Sigma)\le\eta$, then for
every convex set $A$,
\begin{align*}
 \P(W\in A)-\P(Z_{I_d}\in A)
 &\le P:=\eta+(1+M)b/\sqrt8+C_{\scriptscriptstyle\mathrm{GSA}}M\sqrt{\pi b/2},\\
 \P(Z_{I_d}\in A)-\P(W\in A)
 &\le 2\eta+2b/\sqrt8+P+2C_{\scriptscriptstyle\mathrm{GSA}}\sqrt b/\sqrt{2\pi}.
\end{align*}
\end{lemma}
\begin{proof}
The proof first treats a nonempty closed convex set $A$.
The notation is $\mu=\mathcal L(W)$, $Z=Z_{I_d}$, and
\[
 r_A(x)=\dist(B^{-1/2}x,B^{-1/2}A),\qquad x\in\R^d,\qquad
 A_t=\{x\in\R^d:r_A(x)\le t\},\qquad t\ge0.
\]
Since $B^{-1/2}A_t=(B^{-1/2}A)_t$ and $B^{-1/2}Z\sim\mathcal N(0,B^{-1})$,
Lemma~\ref{lem:shell} gives
\[
 \P(Z\in A_t)-\P(Z\in A)\le C_{\scriptscriptstyle\mathrm{GSA}}t\sqrt b,
 \qquad t\ge0.
\]

The convolution density $p_Y(y)=\int_{\R^d}\varphi_B(y-w)\,\mu(dw)$
is positive and continuous, so $p_Y\le M\varphi_\Sigma$ everywhere.
A version of the conditional law given $Y=y$ is $\varphi_B(y-w)\mu(dw)/p_Y(y)$.
With $m=\Sigma^{-1}y$, completing the square and using $y-m=Bm$ gives
\[
 \begin{aligned}
 p_Y(y)\E[e^{v^T(W-m)}\mid Y=y]
 &=e^{v^T(y-m)+v^TBv/2}\,p_Y(y+Bv)\\
 &\le M\varphi_\Sigma(y)e^{v^T(B-B\Sigma^{-1}B)v/2},
 \qquad y,v\in\R^d.
 \end{aligned}
\]
For $m\notin A$, with $x=B^{-1/2}m$ and $\zeta$ the projection of $x$
onto $B^{-1/2}A$, the separating direction $v=B^{-1/2}(\zeta-x)$
satisfies the projection inequality
\[
 v^T(w-m)=(\zeta-x)^T(B^{-1/2}w-\zeta)+\|\zeta-x\|^2
 \ge r_A(m)^2=v^TBv,\qquad w\in A.
\]
Thus, using $B\Sigma^{-1}B\succeq0$,
\[
 p_Y(y)\P(W\in A\mid Y=y)
 \le M\varphi_\Sigma(y)e^{-r_A(m)^2/2},\qquad
 y\in\R^d,\quad m=\Sigma^{-1}y\notin A .
\]

To integrate this estimate, we compare the Gaussian laws with
covariances $\Sigma$ and $\Sigma^{-1}$ to the standard Gaussian.
For each eigenvalue $u$ of $B$, their relative entropy contributions satisfy
\[
 \tfrac12\bigl(u-\log(1+u)\bigr)\le\tfrac14u^2,\qquad
 \tfrac12\bigl(\log(1+u)-u/(1+u)\bigr)\le\tfrac14u^2,\qquad u\ge0.
\]
The derivatives of the left sides are $u/[2(1+u)]$ and
$u/[2(1+u)^2]$, both at most $u/2$. Summing and applying Pinsker's inequality~\cite{pinsker}, \cite[pp.~88--89]{tsybakov} gives
$\tv(Z_\Lambda,Z)\le b/\sqrt8$ for $\Lambda\in\{\Sigma,\Sigma^{-1}\}$.

The function $g_A$ is
$g_A(x)=\Ind_{A^c}(x)e^{-r_A(x)^2/2}$, $x\in\R^d$.
Since $\Sigma A$ is convex and
$\Sigma^{-1}Z_\Sigma\sim\mathcal N(0,\Sigma^{-1})$,
\[
 \P(Y\in\Sigma A)
 \le\eta+\P(Z_\Sigma\in\Sigma A)
 =\eta+\P(Z_{\Sigma^{-1}}\in A).
\]
Integrating the conditional probability against $p_Y$, with the bound
$1$ for $m\in A$ and the conditional bound above for $m\notin A$, gives
\[
 \begin{aligned}
 \P(W\in A)
 &\le \P(Y\in\Sigma A)+M\E g_A(Z_{\Sigma^{-1}})\\
 &\le \P(Z\in A)+\eta+(1+M)b/\sqrt8+M\E g_A(Z).
 \end{aligned}
\]
The second line also uses the total variation bound and $0\le g_A\le1$.
By Tonelli's theorem,
\[
 \E g_A(Z)
 =\int_0^\infty t e^{-t^2/2}\bigl[\P(Z\in A_t)-\P(Z\in A)\bigr]\,dt
 \le C_{\scriptscriptstyle\mathrm{GSA}}\sqrt{\pi b/2}.
\]
This gives $\P(W\in A)-\P(Z\in A)\le P$.

For the reverse inequality, we bound the probability that adding
$Z_B$ moves a point $w\notin A$ into $A$. With $x=B^{-1/2}w$ and
$\zeta$ its projection onto $B^{-1/2}A$, the projection inequality gives
\[
 \{w+Z_B\in A\}\subseteq
 \{(x-\zeta)^TB^{-1/2}Z_B\le-r_A(w)^2\}.
\]
The scalar on the right has variance $r_A(w)^2$ and is centered Gaussian.
Thus $\P(w+Z_B\in A)\le1-\Phi(r_A(w))$. Conditioning on $W$ and applying
Tonelli's theorem to $1-\Phi(r)=\int_r^\infty\varphi_1(t)\,dt$, $r\ge0$, gives
\[
 \begin{aligned}
 \P(Y\in A)
 &\le \mu(A)+\int_0^\infty\varphi_1(t)\mu(A_t\setminus A)\,dt
 =\tfrac12\mu(A)+\int_0^\infty\varphi_1(t)\mu(A_t)\,dt\\
 &\le\tfrac12\bigl[\mu(A)+\P(Z\in A)+P\bigr]
       +C_{\scriptscriptstyle\mathrm{GSA}}\sqrt b/\sqrt{2\pi}.
 \end{aligned}
\]
The last line uses the bound $\P(W\in A_t)-\P(Z\in A_t)\le P$ for the
closed convex sets $A_t$ and the identity
$\int_0^\infty t\varphi_1(t)\,dt=1/\sqrt{2\pi}$ in the enlargement bound.
Since $\P(Y\in A)\ge\P(Z\in A)-\eta-b/\sqrt8$, rearrangement gives the
second inequality.

The empty set is immediate. For any convex set $A$, the upper bound
on $\mu(A)$ follows by applying the closed-set result to $\overline A$,
since $\P(Z\in\overline A)=\P(Z\in A)$ and $\mu(A)\le\mu(\overline A)$.

For the lower bound, the interior $O=\operatorname{int}A$ is approximated
from within. If $O\ne\varnothing$, segments from points of $A$ toward a
point of $O$ show that $\overline O=\overline A$. Thus
$\P(Z\in O)=\P(Z\in A)$. The compact convex sets
\[
 \{x\in\R^d:\|x\|\le j,\ \dist(x,O^c)\ge1/j\},\qquad j\ge1,
\]
increase to $O$. Indeed, the distance condition is closed and says that
the open ball of radius $1/j$ about $x$ is contained in $O$, which is a
convex constraint. The second inequality extends to $O$ by continuity
from below and then to $A$ since $\mu(O)\le\mu(A)$. If $O=\varnothing$, then
$\P(Z\in A)=0$ and the second inequality is immediate.
\end{proof}

\begin{proposition}[One-sided bound]\label{prop:plus}
For $0<\rho\le\rho_*$, every convex set $A$ satisfies
\[
 \P(W\in A)-\P(Z_{I_d}\in A)\le P_\rho,
\]
where
\[
 P_\rho=\eta_\rho
 +\frac{(1+C_{\mathrm{pdf}}(\rho))C_B(\rho)}{\sqrt8}\,\rho^2
 +C_{\scriptscriptstyle\mathrm{GSA}}C_{\mathrm{pdf}}(\rho)
       \sqrt{\frac{\pi C_B(\rho)}{2}}\,\rho .
\]
Here $\eta_\rho$ and $C_B(\rho)$ are given by
Proposition~\ref{prop:smooth}. On this range, $P_\rho<9.36\rho$.
\end{proposition}
\begin{proof}
Proposition~\ref{prop:smooth} gives $Y=W+Z_B$ with
$b=\|B\|_{\mathrm F}\le C_B(\rho)\rho^2<3.211\rho_*^2<1/2$,
$M=C_{\mathrm{pdf}}(\rho)\ge1$, and $\dc(Y,Z_{I_d+B})\le\eta_\rho$.
Lemma~\ref{lem:remove} gives
\[
 \P(W\in A)-\P(Z_{I_d}\in A)
 \le\eta_\rho+\frac{(1+C_{\mathrm{pdf}}(\rho))b}{\sqrt8}
       +C_{\scriptscriptstyle\mathrm{GSA}}C_{\mathrm{pdf}}(\rho)\sqrt{\pi b/2}
 \le P_\rho .
\]
Substitution of $\rho\le\rho_*$ and the bounds from
Proposition~\ref{prop:smooth} and Lemma~\ref{lem:shell}
in the definition of $P_\rho$ gives $P_\rho/\rho<9.36$.
\end{proof}

\begin{proposition}[Convex-set distance for Rademacher sums]\label{prop:frame}
For $\rho>0$, let
\[
 C_{\scriptscriptstyle\mathrm{RS}}(\rho)=
 \begin{cases}
 \displaystyle\min\left\{\frac1\rho,\,
     \max\left\{\frac{4C_{W_2}}{\sqrt\pi}
        \frac{\sqrt{1+\rho}}{(1-\rho)^{3/2}},\ \frac{P_\rho}{\rho}\right\}\right\},
       & 0<\rho\le\rho_*,\\[4pt]
 1/\rho, & \rho>\rho_*,
 \end{cases}
\]
Here $P_\rho$ is given by Proposition~\ref{prop:plus}. For every
deterministic Parseval frame,
\[
 \dc(W,Z_{I_d})\le C_{\scriptscriptstyle\mathrm{RS}}(\rho)\rho .
\]
\end{proposition}
\begin{proof}
For $0<\rho\le\rho_*$, Propositions~\ref{prop:minus} and~\ref{prop:plus}
give the two entries in the maximum. The bound $\dc(W,Z_{I_d})\le1$
supplies the other term $1/\rho$ in the minimum and also applies for
$\rho>\rho_*$.

The function $C_{\scriptscriptstyle\mathrm{RS}}$ is Borel measurable,
as needed when the frame is random in Proposition~\ref{prop:swap}.
Indeed, continuity in $z$ allows the supremum in $\eta_\rho$ to be taken
over rational $z\in(0,1/2]$.
\end{proof}

\section{Bounds for independent summands}\label{sec:independent}
Independent copies of the summands allow us to apply the frame bound
conditionally. For the original sum $S$ and an independent copy $\tilde S$,
the comparison is with $(S+\tilde S)/2+G$ over convex sets, where
$G\sim\mathcal N(0,I_d/2)$ is independent of both. The remaining comparison
with $Z_{I_d}$ is in total variation.

\begin{proposition}[Comparison by random swaps of independent copies]\label{prop:swap}
For independent centered vectors $X_1,\ldots,X_n$ in $\R^d$ with
$\Sigma_i=\Var(X_i)$, $\sum_{i=1}^{n}\Sigma_i=I_d$ and
$\delta^2\coloneqq\sum_{i=1}^{n}\E\|X_i\|^4<\infty$, the sum is
$S\coloneqq\sum_{i=1}^{n} X_i$. Let $(\tilde X_i)_{i=1}^n$ be an independent
copy of the array, $\tilde S=\sum_{i=1}^n\tilde X_i$, and
$G\sim\mathcal{N}(0,I_d/2)$ be independent of both arrays. The parameters
are defined by
\[
 \begin{aligned}
 u&=\delta^{-2}\sum\nolimits_{i=1}^{n}(\tr\Sigma_i)^2,&
 v&=\delta^{-2}\sum\nolimits_{i=1}^{n}\tr(\Sigma_i^2),\\
 \kappa^2&=(1+u+2v)/2,& h^2&=(1+u)/2 ,
 \end{aligned}
\]
where $\kappa,h\ge0$. For $D_i=(X_i-\tilde X_i)/2$, $i=1,\ldots,n$, let
\[
 H=\sum\nolimits_{i=1}^n D_iD_i^T,\qquad
 L=\sum\nolimits_{i=1}^n\|D_i\|^4,\qquad e=\|2H-I_d\| .
\]
On $\{e<1\}$, the matrix $H$ is positive definite and
$\rho_H=(\sum_{i=1}^n\|H^{-1/2}D_i\|^4)^{1/2}$. Then
$4\E L=\kappa^2\delta^2$ and
\[
 \dc\big(S,(S+\tilde{S})/2+G\big)
 \le \E\big[(1-e)\rho_H C_{\scriptscriptstyle\mathrm{RS}}(\rho_H)\,;\ e<1\big]
       +(1+1/\sqrt8)h\delta ,
\]
where the expectation is restricted to $\{e<1\}$ and
$C_{\scriptscriptstyle\mathrm{RS}}$ is defined in Proposition~\ref{prop:frame}.
\end{proposition}
\begin{proof}
For the midpoint $T=(S+\tilde S)/2$ and i.i.d.\ Rademacher variables
$\varepsilon_1,\ldots,\varepsilon_n$ independent of the arrays,
$T+\sum_{i=1}^{n}\varepsilon_iD_i$ has the law of $S$. Indeed, its $i$th
summand is $X_i$ when $\varepsilon_i=1$ and $\tilde X_i$ when
$\varepsilon_i=-1$, and the triples $(X_i,\tilde X_i,\varepsilon_i)$
are independent across $i$.
Conditional on the pairs $(X_i,\tilde X_i)_{i=1}^{n}$, the vectors
$T,D_1,\ldots,D_n$ are fixed, the signs $\varepsilon_i$ remain i.i.d.\
Rademacher variables, and $T+G$ has law $\mathcal N(T,I_d/2)$.

Averaging the conditional bound requires estimates of the fourth-moment
sum $L$ and the covariance deviation $2H-I_d$.
For $i=1,\ldots,n$, the moment quantities are $m_i=\E\|X_i\|^4$,
$t_i=\tr\Sigma_i$ and $s_i=\tr(\Sigma_i^2)$. The total covariance gives $\delta>0$, and
$0\le s_i\le t_i^2\le m_i$ gives $0\le v\le u\le1$. Independence and centering give
\[
 \E D_iD_i^T=\Sigma_i/2,\qquad
 16\E\|D_i\|^4=2m_i+2t_i^2+4s_i.
\]
In the fourth-moment expansion, the terms linear in $X_i^T\tilde X_i$
have zero expectation by conditioning, and
$\E(X_i^T\tilde X_i)^2=\E(X_i^T\Sigma_iX_i)=s_i$.
Since $\|D_iD_i^T\|_{\mathrm F}^2=\|D_i\|^4$ and the matrices
$D_iD_i^T-\Sigma_i/2$ are independent and centered,
\[
 4\E L=\kappa^2\delta^2,\qquad
 \E\|2H-I_d\|_{\mathrm F}^2
 =4\sum\nolimits_{i=1}^{n}\bigl(\E\|D_i\|^4-s_i/4\bigr)=h^2\delta^2.
\]

The comparison is conditional on the pairs. If $e<1$, then
$H\succeq(1-e)I_d/2$ and
\[
 {\textstyle\sum_{i=1}^{n}(H^{-1/2}D_i)(H^{-1/2}D_i)^T}=I_d,\qquad
 (1-e)\rho_H\le 2\sqrt L.
\]
The first identity gives $\rho_H>0$, and the bound on $\rho_H$ follows
from $\|H^{-1/2}\|^2\le2/(1-e)$. By affine invariance and
Proposition~\ref{prop:frame}, the conditional distance to
$\mathcal N(T,H)$ is at most $w=\rho_H C_{\scriptscriptstyle\mathrm{RS}}(\rho_H)$.
With $p=\tv(\mathcal N(0,H),\mathcal N(0,I_d/2))$, the distance to
$\mathcal N(T,I_d/2)$ is therefore at most
\[
 \min\{1,w+p\}\le (1-e)w+(1-e)p+e,
\]
by $\min\{1,y\}\le(1-e)y+e$ for $y\ge0$ and $0\le e<1$.
For the eigenvalues $\lambda_1,\ldots,\lambda_d$ of $2H-I_d$,
\[
 \lambda_j-\log(1+\lambda_j)
 =\lambda_j^2\int_0^1\frac{r}{1+r\lambda_j}\,dr
 \le\frac{\lambda_j^2}{2(1-e)},\qquad j=1,\ldots,d.
\]
Half the sum of the left sides is the relative entropy of
$\mathcal N(0,H)$ with respect to $\mathcal N(0,I_d/2)$.
Pinsker's inequality~\cite[pp.~88--89]{tsybakov} gives
$p\le F/\sqrt{8(1-e)}$, where $F=\|2H-I_d\|_{\mathrm F}$.
With $e\le F$, the conditional distance is at most
$(1-e)w+(1+1/\sqrt8)F$ on $\{e<1\}$. On $\{e\ge1\}$, the distance
is at most $1\le F$. Taking expectations and the supremum over convex
sets gives the stated bound, since $\E F\le h\delta$ by Cauchy--Schwarz.
\end{proof}

\begin{proposition}[Total variation after Gaussian convolution]\label{prop:fixed}
The sum $S\coloneqq\sum_{i=1}^n X_i$ has independent centered summands
$X_i\in\R^d$, with $\sum_{i=1}^{n}\Var(X_i)=I_d$ and
$\delta^2\coloneqq\sum_{i=1}^{n}\E\|X_i\|^4<\infty$.
If $\tilde{S}$ is an independent copy of $S$ and $G\sim\mathcal{N}(0,I_d/2)$
is independent of $(S,\tilde{S})$, then
\[
 \tv\big((S+\tilde{S})/2+G,Z_{I_d}\big)\le a_0\delta+b_0\delta^2,
 \qquad
 a_0={\sqrt3(4-\pi)}/{16},\quad b_0={7\sqrt6(1-\log2)}/{24}.
\]
\end{proposition}
\begin{proof}
uses a fixed Gaussian component in a Lindeberg replacement to quantify the
contribution of higher-order moments to the Gaussian approximation error. We use Slepian inter-
polation and a leave-one-out expansion as in Chernozhukov et al. [10]. The first preprint version
of [11] also uses this interpolation with a fixed Gaussian component.

The proof uses the smooth-case Slepian interpolation with a leave-one-out expansion  
as in Chernozhukov et al.~\cite{cck-maxima,cck-nearly-optimal}. Zhilova~\cite{zhilova-edgeworth} used a fixed Gaussian component with a Lindeberg replacement to quantify the
contribution of higher-order moments to the Gaussian approximation error.

For a Borel set $A$, the smoothed indicator is
$f_s=\Ind_A*\varphi_{sI_d}$, $s>0$.
Lemma~\ref{lem:gaussian-convolution} with $\Sigma=sI_d$ gives
\[
 |\inner{U}{f_s^{(k)}(z)}|\le D_k(s)\|U\|_{\mathrm F},\qquad
 D_k(s)=\frac{\sqrt{k!}}{2s^{k/2}},\qquad z\in\R^d,\quad s>0,\quad k\ge1,
\]
for every order-$k$ tensor $U$.

For an independent copy $(\tilde X_i)_{i=1}^n$ of the array, the
$2n$ independent terms $Y_i=X_i/2$ and $Y_{n+i}=\tilde X_i/2$,
$i=1,\ldots,n$, give
$T=(S+\tilde S)/2=\sum_{j=1}^{2n}Y_j$. With
$\Gamma_j=\E Y_jY_j^T$ and $b_j=\E\|Y_j\|^4$, $j=1,\ldots,2n$,
\[
 \sum\nolimits_{j=1}^{2n}\Gamma_j=I_d/2,\qquad \sum\nolimits_{j=1}^{2n}b_j=\delta^2/8.
\]
The interpolation is $g(t)=\E f_s(V)$, with $s=1-t/2$ and
$V=\sqrt t\,T$, $0\le t\le1$. Since $\Var(T)=I_d/2$, the sum
of $V$ and an independent $\mathcal N(0,sI_d)$ vector has covariance
$(t/2+s)I_d=I_d$.
Thus $g(0)=\P(Z_{I_d}\in A)$ and $g(1)=\P(T+G\in A)$.
The heat equation $\partial_s f_s=\tfrac12\Delta f_s$, where
$\Delta=\sum_{j=1}^d\partial_{jj}$ is the Euclidean Laplacian,
and the chain rule give
\[
 g'(t)=\frac1{2\sqrt t}\sum\nolimits_{j=1}^{2n}\E[Y_j^T\nabla f_s(V)]
       -\frac12\sum\nolimits_{j=1}^{2n}\E\inner{\Gamma_j}{f_s^{(2)}(V)},
 \qquad 0<t<1.
\]
The derivative bounds and $\E\|T\|<\infty$ justify differentiation
under the expectation on compact subintervals of $(0,1)$. Joint continuity of
$f_s(z)$ and $0\le f_s\le1$ give continuity of $g$ on $[0,1]$.

The expansion point $V_j=V-\sqrt t\,Y_j$ is independent of $Y_j$,
and $M_j=\E Y_j^{\otimes3}$, $j=1,\ldots,2n$, is the third-moment tensor.
At $V_j$, a second-order expansion of $Y_j^T\nabla f_s$ and a first-order
expansion of $\inner{\Gamma_j}{f_s^{(2)}}$ isolate the covariance
and third-moment terms.
Independence of $V_j$ and $Y_j$, the identity $\E Y_j=0$, and the
fourth-derivative bound give the remainder bounds
\[
 \begin{aligned}
 &\left|\frac{\E[Y_j^T\nabla f_s(V)]}{2\sqrt t}
  -\frac12\E\inner{\Gamma_j}{f_s^{(2)}(V_j)}
  -\frac{\sqrt t}{4}\E\inner{M_j}{f_s^{(3)}(V_j)}\right|
 \le \frac{t}{12}D_4(s)b_j,\\
 &\frac12\left|\E\inner{\Gamma_j}{f_s^{(2)}(V)-f_s^{(2)}(V_j)}\right|
 \le \frac{t}{4}D_4(s)\|\Gamma_j\|_{\mathrm F}\E\|Y_j\|^2
 \le \frac{t}{4}D_4(s)b_j .
 \end{aligned}
\]
The fourth derivative is contracted with $Y_j^{\otimes4}$ and
$\Gamma_j\otimes Y_j^{\otimes2}$, respectively. The second bound uses
$\|\Gamma_j\|_{\mathrm F}\le\tr\Gamma_j$ and
$(\tr\Gamma_j)^2\le b_j$. Evaluating the third derivatives at the
common point $V$ allows the tensors $M_j$ to be summed. The resulting
error is at most
\[
 \frac{\sqrt t}{4}
 \left|\E\inner{M_j}{f_s^{(3)}(V)-f_s^{(3)}(V_j)}\right|
 \le\frac{t}{4}D_4(s)\|M_j\|_{\mathrm F}\E\|Y_j\|
 \le\frac{t}{4}D_4(s)b_j
\]
by integration along the segment, $\|M_j\|_{\mathrm F}\le\E\|Y_j\|^3$,
and H\"older's inequality.
Summing, the covariance terms cancel and
\[
 \left|g'(t)-\frac{\sqrt t}{4}\inner{M}{\E f_s^{(3)}(V)}\right|
 \le\frac{7t}{96}D_4(s)\delta^2,\qquad
 M={\textstyle\sum_{j=1}^{2n}M_j}=\tfrac14{\textstyle\sum_{i=1}^{n}\E X_i^{\otimes3}}.
\]

For a unit vector $a$ and $k,\ell\in\{1,\ldots,d\}$,
Cauchy--Schwarz and $\sum_{i=1}^{n}\Var(X_i)=I_d$ give
\[
 \left|{\textstyle\sum_{i=1}^{n}\E[(a^TX_i)X_{ik}X_{i\ell}]}\right|^2
 \le {\textstyle\sum_{i=1}^{n}\E X_{ik}^2X_{i\ell}^2} .
\]
The supremum over $a$ followed by the sum over $k,\ell$ gives
$\|4M\|_{\mathrm F}^2\le\delta^2$. Thus
\[
 |g'(t)|\le \frac{\sqrt t}{16}D_3(1-t/2)\delta
        +\frac{7t}{96}D_4(1-t/2)\delta^2,\qquad 0<t<1.
\]
The bound is integrable and $g$ is continuous at the endpoints, so it
bounds $|g(1)-g(0)|$ after integration. Substituting
$u=\sqrt{t/(2-t)}$ and $s=1-t/2$, respectively, gives
\[
 \begin{aligned}
 \int_0^1\sqrt t(1-t/2)^{-3/2}\,dt
 &=4\sqrt2\int_0^1\frac{u^2}{1+u^2}\,du=\sqrt2(4-\pi),\\
 \int_0^1 t(1-t/2)^{-2}\,dt
 &=4\int_{1/2}^1(s^{-2}-s^{-1})\,ds=4(1-\log2).
 \end{aligned}
\]
Using $D_3(s)=\sqrt6/(2s^{3/2})$ and $D_4(s)=\sqrt6/s^2$ yields
$a_0\delta+b_0\delta^2$. The total variation bound follows by taking
the supremum over Borel sets $A$.
\end{proof}

\section{\texorpdfstring{Proof of Theorem~\ref{t}}{Proof of the main theorem}}\label{sec:main-proof}
\begin{proof}
For independent summands, the intermediate vector is
$Y=(S+\tilde S)/2+G$ as in Proposition~\ref{prop:swap}. The triangle inequality and
Propositions~\ref{prop:swap} and~\ref{prop:fixed} give
\[
 \begin{aligned}
 \dc(S,Z_{I_d})
 &\le \dc(S,Y)+\tv(Y,Z_{I_d})\\
 &\le \E\big[(1-e)\rho_H C_{\scriptscriptstyle\mathrm{RS}}(\rho_H)\,;\ e<1\big]
       +\bigl[(1+1/\sqrt8)h+a_0\bigr]\delta+b_0\delta^2.
 \end{aligned}
\]
Here $a_0,b_0$ are the coefficients from Proposition~\ref{prop:fixed}.
For a uniform bound, $(1-e)\rho_H\le2\sqrt L$ on $\{e<1\}$ and
$4\E L=\kappa^2\delta^2$ give
\[
 \E\big[(1-e)\rho_H C_{\scriptscriptstyle\mathrm{RS}}(\rho_H)\,;\ e<1\big]
 \le 2\E\sqrt L\,\sup_{\rho>0}C_{\scriptscriptstyle\mathrm{RS}}(\rho)
 \le \kappa\delta\,\sup_{\rho>0}C_{\scriptscriptstyle\mathrm{RS}}(\rho).
\]
The bounds on $\kappa$ and $h$ below give $A\delta+b_0\delta^2$
in the independent and i.i.d.\ cases, with different values of $A$.
Combining this bound with $\dc\le1$ gives a linear coefficient
through $C=A+b_0/C$. Accordingly, the constants are
\begin{equation}\label{eq:main-constants}
 \begin{aligned}
 A_{\mathrm{ind}}&=\sqrt2\,\sup_{\rho>0}C_{\scriptscriptstyle\mathrm{RS}}(\rho)
                         +1+1/\sqrt8+a_0,\\
 A_{\mathrm{iid}}(d)&=\sqrt{1+1/d}\,\sup_{\rho>0}C_{\scriptscriptstyle\mathrm{RS}}(\rho)
                         +1+1/\sqrt8+a_0,\qquad d\ge2,\\
 C_{\mathrm{ind}}&=\frac{A_{\mathrm{ind}}+\sqrt{A_{\mathrm{ind}}^2+4b_0}}2,\qquad
 C_{\mathrm{iid}}=\frac{A_{\mathrm{iid}}(d)+\sqrt{A_{\mathrm{iid}}(d)^2+4b_0}}2,
 \qquad d\ge2,\\
 C_{\mathrm{iid},1}&=0.9496 .
 \end{aligned}
\end{equation}
Proposition~\ref{prop:frame} gives $C_{\scriptscriptstyle\mathrm{RS}}(\rho)$
in terms of $C_{W_2}$, $\rho_*$, and $P_\rho$. The formula for $P_\rho$
in Proposition~\ref{prop:plus} depends on $C_{\scriptscriptstyle\mathrm{GSA}}$,
$C_{\mathrm{pdf}}(\rho)$, $C_B(\rho)$, and $\eta_\rho$.

Since $0\le v\le u\le1$ by Proposition~\ref{prop:swap}, we have
$\kappa\le\sqrt2$ and $h\le1$. Thus the comparison for independent
summands gives $A_{\mathrm{ind}}\delta+b_0\delta^2$.

For i.i.d.\ isotropic summands with $d\ge2$ and $M_4<\infty$, the
comparison applies to $X_i/\sqrt n$, $i=1,\ldots,n$. In this case,
\[
 \delta^2=M_4/n,\qquad u=d^2/M_4,\qquad v=d/M_4.
\]
Since $M_4\ge(\E\|X_1\|^2)^2=d^2$, we have
$\kappa\le\sqrt{1+1/d}$ and $h\le1$. The comparison gives
$A_{\mathrm{iid}}(d)\delta+b_0\delta^2$.

For $(A,C)=(A_{\mathrm{ind}},C_{\mathrm{ind}})$ or
$(A_{\mathrm{iid}}(d),C_{\mathrm{iid}})$, the definition gives $C=A+b_0/C$. If $\delta\le1/C$,
then $A\delta+b_0\delta^2\le C\delta$. If $\delta>1/C$, then $\dc\le1<C\delta$.

To evaluate the supremum in \eqref{eq:main-constants}, let
\[
 g(\rho)=\frac{4C_{W_2}}{\sqrt\pi}
             \frac{\sqrt{1+\rho}}{(1-\rho)^{3/2}},
 \qquad 0<\rho\le\rho_* .
\]
Both $g(\rho)$ and $\rho g(\rho)$ are strictly increasing. With $C_{W_2}=4$,
the equation $\rho g(\rho)=1$ has a unique solution
$\rho_0\in(1/10.898,1/10)$. Indeed, $\rho g(\rho)$ is continuous and
the endpoint inequalities follow from $3.1415<\pi<22/7$. For
$0<\rho\le\rho_0$, we have $g(\rho)\le1/\rho_0$ and
$P_\rho/\rho<9.36<1/\rho_0$ by Proposition~\ref{prop:plus}.
For $\rho>\rho_0$, the bound $1/\rho<1/\rho_0$ applies. At $\rho=\rho_0$,
the two entries in the minimum in Proposition~\ref{prop:frame} both equal
$1/\rho_0$. The uniform Rademacher coefficient is therefore
\begin{equation}\label{eq:rs-uniform-constant}
 \overline C_{\scriptscriptstyle\mathrm{RS}}
 \coloneqq \sup_{\rho>0}C_{\scriptscriptstyle\mathrm{RS}}(\rho)
 =1/\rho_0<10.898 .
\end{equation}
Proposition~\ref{prop:frame} gives the Rademacher bound with this
coefficient. Its value in \eqref{eq:main-constants} gives
$C_{\mathrm{ind}}<16.88$.

Since the supremum equals $1/\rho_0>0$, $A_{\mathrm{iid}}(d)$ decreases
strictly with $d\ge2$. The function $A\mapsto(A+\sqrt{A^2+4b_0})/2$ is
continuous and strictly increasing for $A>0$. Thus $C_{\mathrm{iid}}$
decreases strictly from its value at $d=2$, which is below $14.81$
by \eqref{eq:main-constants} and \eqref{eq:rs-uniform-constant}.
Its limit is given by the same formula in \eqref{eq:main-constants} with
$A_{\mathrm{iid}}(d)$ replaced by $1/\rho_0+1+1/\sqrt8+a_0$.
Numerically, the value at $d=2$ and the limit are approximately
$14.8075$ and $12.3614$, respectively.

In dimension one, an interval probability is a difference of
distribution-function values or left limits. Thus the convex-set
distance is bounded by twice the
Kolmogorov distance. If $M_4<\infty$, the scalar Berry--Esseen
inequality~\cite{shevtsova} gives
\[
 \dc(S_n,Z_{I_1})
 \le 2\cdot0.4748\, {\E|X_1|^3}/{\sqrt n}
 \le C_{\mathrm{iid},1}\sqrt{M_4/n},
\]
where the last inequality uses $\E|X_1|^2=1$ and Cauchy--Schwarz.
\end{proof}

\section{Proofs of the bounds for general moment orders}\label{sec:finite-moments}

\subsection{Proof of Theorem~\ref{thm:finite-p}}
\label{sec:proof-finite-p}
The proofs of Theorems~\ref{thm:finite-p} and~\ref{thm:moment-constraint}
share the cutoff notation below. The two-block bound in
Lemma~\ref{lem:moment-cutoff} is used for Theorem~\ref{thm:moment-constraint}.
Throughout this subsection, $X_1,\ldots,X_n$ satisfy the independence,
centering and covariance assumptions of Theorem~\ref{t}. With
$\Sigma_i=\Var(X_i)$, the quantities for $t>0$ are
\[
 F=\sum\nolimits_{i=1}^n\E[\|X_i\|^4\Ind_{\{\|X_i\|\le t\}}],\qquad
 T_k=\sum\nolimits_{i=1}^n\E[\|X_i\|^k\Ind_{\{\|X_i\|>t\}}],\qquad
 D=\sum\nolimits_{i=1}^n\|\Sigma_i\|_{\mathrm F}^2.
\]

\begin{lemma}[A bound with a cutoff]\label{lem:moment-cutoff}
If $0<t\le1$, $T_2\le1/6$, and $F\le1/3600$, then
\begin{equation}\label{eq:moment-cutoff}
 \dc\!\left({\textstyle \sum_{i=1}^n X_i},Z_{I_d}\right)
 \le \frac{5C_{\mathrm{ind}}}{\sqrt3}\sqrt F
 +\left(225+\left({40}/{3}\right)^2\right)F
 +\frac{10\sqrt6}{3}T_3+\sqrt{24}(F+D).
\end{equation}
\end{lemma}
\begin{proof}
The bound is immediate if $T_3=\infty$.

The cutoff indicator is $J_i=\Ind_{\{\|X_i\|\le t\}}$, with
$p_i=\E J_i$. The mean and covariance of $X_i$ conditional on $J_i=1$
are $\mu_i$ and $C_i$, respectively, both zero if $p_i=0$. The Gaussian
vectors $G_i\sim\mathcal N(0,C_i)$ are independent of one another and
of all $X_i$. The corresponding variables and covariance matrices are
\[
 Y_i=J_i(\mu_i+G_i)+(1-J_i)X_i,\qquad
 Q=\sum\nolimits_{i=1}^n J_i C_i,\qquad \bar Q=\E Q=\sum\nolimits_{i=1}^n p_i C_i.
\]
The $Y_i$ are independent and match the means and covariances of the $X_i$.
The conditional covariance identity gives
\[
 I_d-\bar Q=\sum\nolimits_{i=1}^n\bigl(\E[X_iX_i^T(1-J_i)]+p_i\mu_i\mu_i^T\bigr)
 \succeq0.
\]
For $p_i\ge1/2$, centering and Cauchy--Schwarz give
$p_i\|\mu_i\|^2\le (1-p_i)p_i^{-1}\E[\|X_i\|^2(1-J_i)]$.
For $0<p_i<1/2$, we have $p_i\|\mu_i\|^2\le p_it^2
\le \E[\|X_i\|^2(1-J_i)]$. Hence
$\tr(I_d-\bar Q)\le2T_2$ and $\bar Q\succeq(2/3)I_d$.
Also, $\|C_i\|_{\mathrm F}\le\tr C_i\le\E[\|X_i\|^2\mid J_i=1]$
when $p_i>0$. Jensen's inequality gives
$\sum_{i=1}^n p_i\|C_i\|_{\mathrm F}^2\le F$.
For each index set $s\subseteq\{1,\ldots,n\}$, the covariance in that
block is $Q_s=\sum_{i\in s}J_iC_i$. Independence gives
\[
 \E\|Q_s-\E Q_s\|_{\mathrm F}^2
 =\sum\nolimits_{i\in s}p_i(1-p_i)\|C_i\|_{\mathrm F}^2\le F.
\]

In a random partition $s_1,s_2$, each index is assigned independently
to either block with probability $1/2$. The expected squared deviation
of $\sum_{i\in s_1}p_iC_i$ from $\bar Q/2$ in Frobenius norm is
$\tfrac14\sum_{i=1}^n p_i^2\|C_i\|_{\mathrm F}^2\le F/4$.
Consequently, there is a partition with deviation at most $\sqrt F/2$.
For this partition the second block has the same deviation. Since
$\sqrt F/2\le1/120$, both $\E Q_{s_j}\succeq(13/40)I_d$. Markov's inequality gives
\begin{equation}\label{eq:moment-cutoff-bad}
 \P(Q\not\succeq (3/5)I_d)\le225F,\qquad
 \P(Q_{s_j}\not\succeq I_d/4)\le(40/3)^2F,\qquad j=1,2.
\end{equation}

Conditional on all $(J_i,(1-J_i)X_i)$, $1\le i\le n$, the terms with
$J_i=1$ are independent, each with its law conditional on $J_i=1$.
If $Q\succeq(3/5)I_d$, Theorem~\ref{t} applies to these terms after
centering by $\mu_i$ and multiplication by $Q^{-1/2}$. For any vector
$U$ with finite fourth moment, the centered vector is $V=U-\E U$.
The nonnegative parameters $a,b$ satisfy
$a^2=\E\|V\|^4$ and $b^2=\E\inner{V}{\E U}^2$.
Expansion and Cauchy--Schwarz, with
$b^2\le\|\E U\|^2\E\|V\|^2$, give
\[
 \E\|U\|^4\ge a^2-4ab+6b^2
 =\frac{a^2}{3}+6(b-a/3)^2\ge\frac{a^2}{3}.
\]
Applying this to the conditional law gives
$\E[\|X_i-\mu_i\|^4\mid J_i=1]\le3\E[\|X_i\|^4\mid J_i=1]$.
The square root of the fourth-moment sum after normalization is thus
at most
\[
 \frac{5}{\sqrt3}\left({\textstyle\sum_{i=1}^n J_i\E[\|X_i\|^4\mid J_i=1]}\right)^{1/2},
\]
where the product is zero if $p_i=0$. The conditional Gaussian
comparison law, after undoing the affine transformation and restoring
the fixed tail terms, is the conditional law of $\sum_{i=1}^nY_i$.
Jensen's inequality bounds the average, and
\eqref{eq:moment-cutoff-bad} bounds the probability of the exceptional
event. Together they give
\[
 \dc\!\left({\textstyle \sum_{i=1}^n X_i},{\textstyle \sum_{i=1}^n Y_i}\right)
 \le \frac{5C_{\mathrm{ind}}}{\sqrt3}\sqrt F+225F.
\]

During replacements in either block, the other block supplies the
Gaussian convolution. The comparison of the modified sum with the
target Gaussian uses $f=\Ind_A*\varphi_{Q_0}$ for a convex set $A$
and $Q_0\succeq I_d/4$.
Lemma~\ref{lem:gaussian-convolution} gives
$\sup_{x\in\R^d}\|f^{(k)}(x)\|_{\mathrm F}\le L_k:=2^{k-1}\sqrt{k!}$
for $1\le k\le4$.

For $H\succeq0$ and $u>0$, Gaussian integration by parts on the range
of $H$ gives the first and second derivatives of $\E f(y+Z_{uH})$
by applying $\tfrac12\sum_{k,\ell=1}^dH_{k\ell}\partial_{k\ell}$ once
and twice, respectively, to $f$ inside the expectation.
Bounded continuous derivatives give their limits at $u=0$.
Taylor's formula therefore gives
\[
 \left|\E f(y+Z_H)-f(y)-\tfrac12\inner{H}{f^{(2)}(y)}\right|
 \le \tfrac18L_4\|H\|_{\mathrm F}^2,\qquad y\in\R^d.
\]
With $A_i=J_i\mu_i+(1-J_i)X_i$ and $H_i=J_iC_i$, the law of
$Y_i$ conditional on $(A_i,J_i)$ is Gaussian with mean $A_i$ and
covariance $H_i$. The comparison uses the preceding bound at
$y=x+A_i$, the second-order expansion of $f(x+A_i)$ at $x$, and
$\|f^{(2)}(x+A_i)-f^{(2)}(x)\|_{\mathrm F}\le L_3\|A_i\|$.
Because $\E A_i=0$ and $\E(A_iA_i^T+H_i)=\Sigma_i$, the terms through
order two agree with those for $Z_{\Sigma_i}$. The same bound applied
to $Z_{\Sigma_i}$ gives
\[
 \begin{aligned}
 |\E f(x+Y_i)-\E f(x+Z_{\Sigma_i})|
 &\le \tfrac16L_3\E\|A_i\|^3
   +\tfrac12L_3\E[\|A_i\|\|H_i\|_{\mathrm F}]\\
 &\quad+\tfrac18L_4\bigl(\E\|H_i\|_{\mathrm F}^2+\|\Sigma_i\|_{\mathrm F}^2\bigr),
 \qquad x\in\R^d .
 \end{aligned}
\]
Since $\|\mu_i\|\le t$, $\|C_i\|_{\mathrm F}\le t^2$, and
$p_i\|\mu_i\|\le\E[\|X_i\|(1-J_i)]$, the first two expectations sum
over $i$ to at most $T_3+t^2T_1\le2T_3$ and $t^2T_1\le T_3$,
respectively. Also $\sum_{i=1}^n\E\|H_i\|_{\mathrm F}^2\le F$.
The sum of the replacement bounds is
$(5L_3/6)T_3+(L_4/8)(F+D)$.

The first stage replaces $Y_i$ by independent $Z_{\Sigma_i}$ for
$i\in s_1$, one at a time. Conditional on $(J_i,(1-J_i)X_i)_{i\in s_2}$, the
second block is Gaussian with covariance $Q_{s_2}$ and a fixed mean,
independently of the first block. If $Q_{s_2}\succeq I_d/4$, it supplies
the convolution in the replacement bound. On the complementary event, the conditional error is at most one. The new Gaussian block
$s_1$ then supplies covariance
$\sum_{i\in s_1}\Sigma_i\succeq\E Q_{s_1}\succeq(13/40)I_d$ while
the terms in $s_2$ are replaced. The final sum has law $\mathcal N(0,I_d)$.
The two block errors, the exceptional probability from
\eqref{eq:moment-cutoff-bad}, and the error from replacing the terms
at or below the cutoff together give \eqref{eq:moment-cutoff}.
\end{proof}

\begin{proof}[Proof of Theorem~\ref{thm:finite-p}]
Let $S=\sum_{i=1}^n X_i$, $\eta=B_q^{1/(q-2)}>0$,
$r_q=3^{(3-q)/(q-2)}$, and $t=3^{1/(q-2)}\eta$.
The fraction of $B_q$ contributed by the terms at or below $t$ is
\[
 s=B_q^{-1}\sum_{i=1}^n\E[\|X_i\|^q\Ind_{\{\|X_i\|\le t\}}].
\]
Bounding powers of $\|X_i\|$ separately on the two regions gives
\[
 T_2\le(1-s)/3,\qquad F\le3s\,r_q^2\eta^2,\qquad
 T_3\le(1-s)r_q\eta .
\]
In the notation of the proof of Lemma~\ref{lem:moment-cutoff}, the
cutoff data are $(J_i,(1-J_i)X_i)_{i=1}^n$. The conditional means of
the $X_i$ given these data are $A_i=J_i\mu_i+(1-J_i)X_i$.
The $A_i$ are independent and centered, and the conditional covariance
of $S$ is $Q=\sum_{i=1}^n J_i C_i$. By total covariance,
$\sum_{i=1}^n\Var(A_i)=I_d-\bar Q$, where $\bar Q=\E Q$.
The estimates in that proof hold for every cutoff $t>0$ and give
\[
 \bar Q\succeq (1-2T_2)I_d,\qquad
 \E\|Q-\bar Q\|_{\mathrm F}^2\le F.
\]
For the chosen cutoff, $T_2\le(1-s)/3$, so
$\bar Q\succeq mI_d$ with $m=(1+2s)/3$.

The centered fourth moment of a vector $V$ with finite fourth moment
is at most $c\,\E\|V\|^4$, where $c=1+2/\sqrt3$.
For a scalar $V$ with $\E V\ne0$, the substitution $x=V/\E V$ in the
identity
\[
 c x^4-(x-1)^4-(10+6\sqrt3)(x-1)
 =\frac2{\sqrt3}\left(x^2+\sqrt3x-\frac{3(1+\sqrt3)}2\right)^2,
 \qquad x\in\R,
\]
gives the inequality after taking expectations and multiplying by
$(\E V)^4$, since $\E(x-1)=0$. For mean zero it holds since $c>1$.
The vector case follows by applying the scalar inequality to $g^TV$ for each
$g\in\R^d$ and integrating with respect to $g\sim\mathcal N(0,I_d)$,
using $\E_g|g^Tv|^4=3\|v\|^4$ for $v\in\R^d$.
Consequently,
\[
 U=\sum_{i=1}^n J_i\,\E[\|X_i-\mu_i\|^4\mid J_i=1],
 \qquad \E U\le c F ,
\]
with a zero summand if $\P(J_i=1)=0$.

Conditional on the cutoff data, the retained summands are independent
with the laws of the corresponding $X_i$ given $J_i=1$. Let
$K=\bar Q^{-1/2}(Q-\bar Q)\bar Q^{-1/2}$ and $e=\|K\|$.
For $e<1$, $Q\succeq m(1-e)I_d$. Theorem~\ref{t}, applied
after centering and normalization by $Q^{-1/2}$, gives conditional
convex-set distance at most $a\sqrt U/[m(1-e)]$ to the Gaussian
with mean $\sum_{i=1}^n A_i$ and covariance $Q$. Here $a=C_{\mathrm{ind}}$.
For i.i.d.\ summands the retained conditional laws are identical.
Their number $N$ is positive since $Q\succ0$. The normalized conditional
sum is $N^{-1/2}$ times the sum of the retained vectors centered by their
common mean and multiplied by $\sqrt N\,Q^{-1/2}$. These transformed
vectors are isotropic and i.i.d., so the same bound holds with
$a=C_{\mathrm{iid}}$ for $d\ge2$ and $a=C_{\mathrm{iid},1}$ for $d=1$.

After normalization by $\bar Q^{-1/2}$, the Gaussian covariance
comparison in the proof of Proposition~\ref{prop:swap} gives total
variation distance at most $\|K\|_{\mathrm F}/\sqrt{8(1-e)}$
between the Gaussians with covariances $Q$ and $\bar Q$ and the same
mean. The inequality $\min\{1,y\}\le(1-e)y+e$ for $y\ge0$ and
$0\le e<1$ therefore bounds the conditional error to the Gaussian
with covariance $\bar Q$ by
\[
 \frac{a\sqrt U}{m}+(1+1/\sqrt8)\|K\|_{\mathrm F}.
\]
On $\{e\ge1\}$, this bound holds since $\|K\|_{\mathrm F}\ge1$.
Averaging the conditional bound, using
$\E\|K\|_{\mathrm F}\le\sqrt F/m$ and $\E\sqrt U\le\sqrt{cF}$,
and taking the supremum over convex sets gives
\[
 \dc\!\left(S,{\textstyle\sum_{i=1}^n A_i}+Z_{\bar Q}\right)
 \le P(a)\frac{\sqrt F}{m},\qquad
 P(a)=\sqrt c\,a+1+1/\sqrt8,\qquad a>0 ,
\]
where $Z_{\bar Q}$ is independent of the $A_i$. This Gaussian
smooths the indicators for the remaining comparison with $Z_{I_d}$.

For a convex set $A$, the function $f=\Ind_A*\varphi_{\bar Q}$
satisfies $\sup_{x\in\R^d}\|f^{(3)}(x)\|_{\mathrm F}
\le\sqrt6/(2m^{3/2})$ by Lemma~\ref{lem:gaussian-convolution}
and $\bar Q\succeq mI_d$. The vectors
$G_i'\sim\mathcal N(0,\Var(A_i))$, $i=1,\ldots,n$, are mutually
independent and independent of $(A_1,\ldots,A_n,Z_{\bar Q})$.
Then $\sum_{i=1}^n G_i'+Z_{\bar Q}\sim\mathcal N(0,I_d)$.
Replacing the $A_i$ by the $G_i'$ one
at a time and expanding $f$ through order two gives
\[
 \dc\!\left({\textstyle\sum_{i=1}^n A_i}+Z_{\bar Q},Z_{I_d}\right)
 \le\frac{\sqrt6}{12m^{3/2}}
       \sum_{i=1}^n(\E\|A_i\|^3+\E\|G_i'\|^3).
\]
The terms through order two cancel because the means and covariances
match. Cauchy--Schwarz and the Gaussian fourth-moment formula give
$\E\|G_i'\|^3\le\sqrt3\,(\E\|A_i\|^2)^{3/2}
\le\sqrt3\,\E\|A_i\|^3$. Since
$\P(J_i=1)\|\mu_i\|\le\E[\|X_i\|(1-J_i)]$ and $\|\mu_i\|\le t$,
we have $\sum_{i=1}^n\E\|A_i\|^3\le2T_3$. Together the two comparisons give
\[
 \dc(S,Z_{I_d})\le P(a)\frac{\sqrt F}{m}
     +B\frac{T_3}{m^{3/2}},\qquad B=\frac{\sqrt6(1+\sqrt3)}6.
\]

With $x=1/m\in[1,3]$, the moment bounds at the chosen cutoff give
\[
 \frac{\dc(S,Z_{I_d})}{r_q\eta}
 \le P(a)\sqrt{\frac32x(3-x)}+\frac{3B}{2}(x-1)\sqrt x.
\]
Since $x(3-x)=9/4-(x-3/2)^2$ and $\sqrt x\le\sqrt3$, Cauchy--Schwarz
bounds the right-hand side by
$(3/2)\sqrt{(3/2)P(a)^2+b^2}+b/2$, where $b=3\sqrt3\,B/2$.
Thus the coefficients are
\begin{equation}\label{eq:finite-moment-constant}
 \begin{aligned}
 \mathcal K_q(a)
 &=3^{(3-q)/(q-2)}
   \left[\frac32\sqrt{\frac32 P(a)^2+b^2}+\frac b2\right],
    \qquad a>0,\quad 3\le q<4,\\
 C_q^{\mathrm{ind}}&=\mathcal K_q(C_{\mathrm{ind}}),\qquad
 C_q^{\mathrm{iid}}=
 \begin{cases}
 \mathcal K_q(C_{\mathrm{iid},1}),&d=1,\\
 \mathcal K_q(C_{\mathrm{iid}}),&d\ge2 .
 \end{cases}
 \end{aligned}
\end{equation}
The expression increases with $a$. We have $b<2.90$,
$P(16.88)<26.14$, and $P(14.81)<23.12$, which give the bounds
$C_q^{\mathrm{ind}}<49.7r_q$ and $C_q^{\mathrm{iid}}<44.2r_q$.
For the i.i.d.\ bound, the comparison is applied to $X_i/\sqrt n$,
so $B_q=n^{1-q/2}M_q$ and $B_q^{1/(q-2)}=M_q^{1/(q-2)}/\sqrt n$.

The factor $r_q=3^{1/(q-2)-1}$ decreases strictly from $1$ at $q=3$
to $1/\sqrt3$ as $q\uparrow4$. Thus, for fixed $d$, each coefficient
has its maximum at $q=3$ and approaches its infimum as $q\uparrow4$.
Since $\mathcal K_q$ is strictly increasing in $a$ and $C_{\mathrm{iid}}$
decreases strictly with $d\ge2$ by Theorem~\ref{t}, the maximum of
$C_q^{\mathrm{iid}}$ over $d\ge2$ and $3\le q<4$ occurs at $(d,q)=(2,3)$.
The infimum is approached as $d\to\infty$ and $q\uparrow4$.
The values obtained from \eqref{eq:main-constants} and
\eqref{eq:finite-moment-constant}, rounded to four decimal places, are
\[
 \begin{array}{lrr}
 &\text{Infimum}&\text{Maximum}\\
 C_q^{\mathrm{ind}}&28.6517&49.6262\\
 C_q^{\mathrm{iid}},\ d=1&4.6823&8.1100\\
 C_q^{\mathrm{iid}},\ d=2&25.4548&44.0889\\
 C_q^{\mathrm{iid}},\ \text{across }d\ge2&21.6697&44.0889
 \end{array}
\]
\end{proof}

\subsection{Proof of the upper bounds in Theorem~\ref{thm:moment-constraint}}
\label{sec:proof-moment-upper}
\begin{proof}
For $3\le q<4$, the first estimate bounds the $q$th moment above
$2\sqrt d$ in terms of $s=M_q-d^{q/2}$. With $\alpha=q/2$,
for $x\ge4$ and $\alpha\in[3/2,2]$,
\[
 x^\alpha-1-\alpha(x-1)\ge (5/16)x^\alpha .
\]
Indeed, the ratio on the left to $x^\alpha$ is increasing in $x\ge1$,
with derivative $\alpha(\alpha-1)(x-1)x^{-\alpha-1}$. At $x=4$ the
ratio is $1-4^{-\alpha}(1+3\alpha)$. Its derivative in $\alpha$ is
$4^{-\alpha}((1+3\alpha)\log 4-3)>0$ on $[3/2,2]$, and its value at
$\alpha=3/2$ is $5/16$.
The left side is nonnegative for every $x\ge0$ by convexity. Substituting
$x=\|X_1\|^2/d$ and using $\E x=1$ gives
\[
 \E[\|X_1\|^q\Ind_{\{\|X_1\|>2\sqrt d\}}]\le (16/5)s .
\]

The two terms in the bound are $x_0=d/\sqrt n$ and
$y_0=s^{1/(q-2)}/\sqrt n$, with $R=x_0+y_0$.
If $R>1/1024$, the result follows from $\dc\le1$. Otherwise the cutoff
for $X_i/\sqrt n$ in Lemma~\ref{lem:moment-cutoff} is
$t=20R\ge2\sqrt{d/n}$. The last display gives
$T_k\le (16/5)t^{k-q}y_0^{q-2}$ for $k=2,3$.
The contribution to $F$ from $\|X_i/\sqrt n\|\le2\sqrt{d/n}$
is at most $4(d/n)\sum_{i=1}^n\E\|X_i/\sqrt n\|^2=4x_0^2$.
For the contribution from $2\sqrt{d/n}<\|X_i/\sqrt n\|\le t$,
the tail estimate gives $(16/5)t^{4-q}y_0^{q-2}$. Thus
\[
 \begin{aligned}
 F&\le4x_0^2+(16/5)t^{4-q}y_0^{q-2}
       \le 4x_0^2+64R y_0\le64R^2,\\
 T_2&\le4/25,\qquad T_3\le(16/5)R,\qquad D=d/n\le R^2 .
 \end{aligned}
\]
In the first line, $y_0\le R$ and $q-2\ge1$ give
$t^{4-q}y_0^{q-2}\le20R y_0$. These bounds meet the cutoff lemma's
hypotheses. Substitution in \eqref{eq:moment-cutoff}, with
$C_{\mathrm{ind}}<16.88$ and $R\le1/1024$, bounds the error by $1024R$.
Since $M_q\le \widetilde{M}_q$, this proves the upper bound with
$d+(\widetilde{M}_q-d^{q/2})^{1/(q-2)}$.

At $q=4$, Theorem~\ref{t} applies directly, with
$\sqrt{\widetilde{M}_4}\le d+\sqrt{\widetilde{M}_4-d^2}\le\sqrt{2\widetilde{M}_4}$. For $q>4$, H\"older's inequality
with $\E\|X_1\|^2=d$ gives
\begin{equation}\label{eq:moment-interpolation}
 \sqrt{M_4}\le d^{(q-4)/(2(q-2))}M_q^{1/(q-2)} .
\end{equation}
Using $M_q\le\widetilde{M}_q$ and $\dc\le1$, Theorem~\ref{t} gives the claimed
upper bound.
The same interpolation with respect to the finite measure
$\sum_{i=1}^n\mathcal L(X_i)$ gives
$\delta\le d^{(q-4)/(2(q-2))}B_q^{1/(q-2)}$ for independent
summands with total covariance $I_d$.
\end{proof}

\appendix

\section{Proofs of lower bounds}\label{sec:lower-bounds}
In this appendix $c,C$ denote positive universal constants that may
change from line to line.
\subsection{Lattice lower bounds}
\label{sec:lattice-projection}
The cosine-basis law below gives the lattice bound used in
Theorem~\ref{thm:moment-constraint}. The permutation construction
gives the coordinate-moment bound in Theorem~\ref{thm:simplex-lower}(ii).

If $\E T^2=1$ and $T$ is supported on a lattice with spacing $h>0$,
then $\P(|T|\le2)\ge3/4$. The interval $[-2,2]$ contains at most
$1+4/h$ lattice points, so some $t$ satisfies
\begin{equation}\label{eq:lattice-atom}
 \P(T=t)\ge\frac{3h}{4(h+4)}\ge\frac3{20}\min\{1,h\}.
\end{equation}
If $T$ is a unit projection of $S_n$, the corresponding affine
hyperplane has Gaussian mass zero, so \eqref{eq:lattice-atom} gives
a lower bound on $\dc(S_n,Z_{I_d})$.

For $d\ge2$, the vectors $u_1,\ldots,u_d$ are the columns of the
orthogonal matrix $O$ given by
\[
 O_{k1}=d^{-1/2},\qquad
 O_{kj}=\sqrt{2/d}\cos\!\bigl(\pi(j-1)(k-\tfrac12)/d\bigr),
 \quad 2\le j\le d,\quad 1\le k\le d.
\]
For $t=\ell\pi/d$ with $1\le\ell<2d$, the sum
$\sum_{k=1}^d\cos((k-\tfrac12)t)=\sin(dt)/(2\sin(t/2))$ vanishes.
The product-to-sum formula then gives $O^TO=I_d$.
For $L\ge\sqrt d$, the scale is $b^2=(d-1)/(1-L^{-2})\le d$.
The law of $X$ assigns probability $1/(2L^2)$ to each of $\pm Lu_1$
and probability $1/(2b^2)$ to each of $\pm bu_j$, $2\le j\le d$.
The probabilities sum to one. Symmetry gives $\E X=0$ and each pair
contributes $u_j u_j^T$ to $\E XX^T=I_d$. For $d=1$, the law assigns
probability $1/(2L^2)$ to each of $\pm L$ and the remaining probability
to zero.
For every $q\ge3$, the law satisfies
\begin{equation}\label{eq:lattice-moments}
 \begin{gathered}
 \|X\|\le L\ \text{almost surely},\qquad
 \E\|X\|^q\le d^{q/2}+L^{q-2},\\
 \max_{1\le k\le d}\E|X_k|^q\le L^{q-2}d^{-q/2}+2^{(q-2)/2}.
 \end{gathered}
\end{equation}
Indeed, for $d\ge2$ we have $M_q=L^{q-2}+(d-1)b^{q-2}$. The contribution to the coordinate moment
from $j\ge2$ is
$b^{q-2}\sum_{j=2}^d|O_{kj}|^q
 \le b^{q-2}(2/d)^{(q-2)/2}\sum_{j=2}^d O_{kj}^2
 \le 2^{(q-2)/2}$.
The dimension-one assertions follow from the masses directly.
The projection on $u_1$ (on $1$ if $d=1$) of the normalized sum
lies in $(L/\sqrt n)\mathbf Z$. By \eqref{eq:lattice-atom},
\begin{equation}\label{eq:lattice-lower}
 \dc(S_n,Z_{I_d})\ge \frac3{20}\min\{1,L/\sqrt n\}.
\end{equation}

\begin{proof}[Proof of Theorem~\ref{thm:simplex-lower}(ii)]
The scale is $L=(\widetilde{M}_q-d^{q/2})^{1/(q-2)}\ge\sqrt d$
and the direction is $u=d^{-1/2}(1,\ldots,1)^T$. The law of $X$ has
mass $1/(2L^2)$ at each of $Lu$ and $-Lu$. For $d=1$, the remaining
mass is at zero.

For $d\ge2$, $V$ is a uniformly random permutation of a zero-sum
vector of squared norm $d-1$, with entries bounded by $1$.
In even dimension the vector consists of equal numbers of
$\sqrt{(d-1)/d}$ and $-\sqrt{(d-1)/d}$. In odd dimension it consists
of equal numbers of $1$ and $-1$, and one zero. Exchangeability and
the zero sum give $\E V=0$, $\E V_j^2=(d-1)/d$, and
$\E V_jV_k=-1/d$ for $j\ne k$, so $\Var(V)=I_d-uu^T$.
The remaining component of the law of $X$ has weight $1-L^{-2}$
and is the law of $V/\sqrt{1-L^{-2}}$.

The two parts contribute $uu^T$ and $I_d-uu^T$ to the covariance,
respectively, so $X$ is centered and isotropic. The second part
has squared norm $(d-1)/(1-L^{-2})\le d$, which gives
$\|X\|\le L$ and $\E\|X\|^q\le L^{q-2}+d^{q/2}=\widetilde{M}_q$.
For $d\ge2$, the contribution of the second part to
$\E|X_j|^q$ is at most
\[
 (1-L^{-2})^{1-q/2}\E V_j^2
 \le \left(\frac{d-1}{d}\right)^{2-q/2}\le 1 ,
\]
since $|V_j|\le1$. Adding the contribution of the two atoms gives
$\E|X_j|^q\le 1+L^{q-2}d^{-q/2}=\widetilde{M}_q/d^{q/2}$.
For $d=1$, the two atoms give the entire moment.

The projection $u^TS_n$ has variance one and lies on a lattice of
spacing $L/\sqrt n$. By \eqref{eq:lattice-atom}, some affine
hyperplane has $S_n$-mass at least $\frac3{20}\min\{1,L/\sqrt n\}$
and Gaussian mass zero. Finally,
$L\ge 2^{-1/(q-2)}\widetilde{M}_q^{1/(q-2)}
\ge \frac12\widetilde{M}_q^{1/(q-2)}$, which proves the claimed lower bound.
\end{proof}

\subsection{Proof of Theorem~\ref{thm:simplex-lower}(i)}
\label{sec:proof-simplex}
\begin{proof}
For $d=1$, the law assigns mass $1/(2a^2)$ to each of $\pm a$ and
the remaining mass to zero. This gives the stated moment. By \eqref{eq:lattice-atom} with
spacing $a/\sqrt n$, a singleton has the required $S_n$-mass and has
Gaussian mass zero.

The construction for $d\ge2$ uses $N=d+1\ge3$,
$H=\mathbf 1^\perp\subset\R^N$, $P=I_N-N^{-1}\mathbf 1\mathbf 1^T$,
and $u_j=Pe_j$, $j=1,\ldots,N$.
The vector $V$ is uniform on $\{\sqrt N\,u_j:1\le j\le N\}$. Since
$\sum_{j=1}^N u_j=0$ and $\sum_{j=1}^N u_ju_j^T=P$, we have $\E V=0$,
$\Var(V)=P$, and $\|V\|^2=d$. The vector $X$ equals $aV$ with
probability $a^{-2}$ and zero otherwise, with the choice independent
of $V$. Under an orthonormal identification of $H$ with $\R^d$, $X$
has the required covariance and moment. For $a=1$, the basis of $H$
consists of the last $d$ columns of the $N\times N$ cosine basis from
Section~\ref{sec:lattice-projection}.
The coordinates of $\sqrt N\,u_j$ are $\sqrt N$ times the corresponding
entries in row $j$ of that basis, so their absolute values are at most
$\sqrt2$.

We work on $H$ and compare with $Z_P$. The isometry preserves the
convex-set distance. The projection of $S_n$ on $(e_1-e_2)/\sqrt2$
has variance one and lies on a lattice with spacing
$\varepsilon/\sqrt2$, where $\varepsilon=a\sqrt{N/n}$. By
\eqref{eq:lattice-atom}, a hyperplane in $H$ has $S_n$-mass at least
$c\min\{1,\varepsilon\}$ and $Z_P$-mass zero. This proves the lower
bound when $\varepsilon\ge\varepsilon_0$, for any fixed
$\varepsilon_0\in(0,1)$, with a constant depending on $\varepsilon_0$.

For $\varepsilon<\varepsilon_0$, the lower bound comes from the cubic
Taylor term of a smooth test with convex sublevel sets. For a standard
normal $\xi$, the functions are
\[
 f(x)=\E(x+\xi-1)_+,\quad x\in\R,\qquad
 \chi(x)=\sum\nolimits_{j=1}^{N}f(x_j),\quad x\in H.
\]
Then $f'=\Phi(\,\cdot-1)$ and $f''=\varphi_1(\,\cdot-1)$, so $f$ is
convex and $1$-Lipschitz, with all positive-order derivatives bounded.
For a universal constant $L_0\ge1$ specified below, the test function is
\[
 \mu=\E \chi(Z_P),\qquad b=L_0(\sqrt N+N\varepsilon),\qquad
 h(x)=\Phi((\chi(x)-\mu)/b),\quad x\in H.
\]
The sublevel sets of $h$ are convex, and integration over levels
$t\in(0,1)$ gives
\[
 |\E h(S_n)-\E h(Z_P)|\le \dc(S_n,Z_P).
\]

The difference of expectations is estimated by replacing the
normalized copies of $X$ by independent copies of $Z_P/\sqrt n$ one
at a time. In each comparison, $U$ is the sum of the other $n-1$
terms and is independent of both candidate terms.
The law of $U$ is invariant under coordinate permutations.
The scalar Gaussian $\xi_\sigma$ has variance $\sigma^2=P_{11}=1-1/N$
and mean zero.
For $g=f$ or $g=f^{(3)}$,
\begin{equation}\label{eq:simplex-scalar-replacement}
 |\E g(U_1)-\E g(\xi_\sigma)|\le C\varepsilon.
\end{equation}
Indeed, a scalar replacement in coordinate $1$ cancels the terms
through degree two and costs at most $\|g^{(3)}\|_\infty/6$ times
the sum of the third absolute moments of the two normalized
terms. Since $|V_1|\le\sqrt N$ and $\E V_1^2=\sigma^2\le1$,
$\E|e_1^T X|^3\le a\sqrt N$. Thus at most $n-1$ scalar replacements cost
$C a\sqrt{N/n}$. Adding the omitted centered Gaussian costs
at most $\|g''\|_\infty/(2n)\le C\varepsilon$.

The Efron--Stein inequality gives $\Var(\chi(U))\le CN$.
Indeed, $\chi$ is $1$-Lipschitz in the $\ell_1$ norm and
$\E\|X\|_1^2=\E\|V\|_1^2\le4N$ because $\|u_j\|_1\le2$.
Resampling one original summand therefore gives a squared change
in $\chi$ with expectation at most $CN/n$.
The independent Gaussian part of $U$ can be represented as
$sP Z_{I_N}$ with $0\le s\le1$. Resampling the $j$th coordinate of $Z_{I_N}$ changes
$U$ in direction $u_j$ and gives an expected squared change in $\chi$
bounded by a universal constant.
Summing the resampling bounds proves the variance estimate.
By \eqref{eq:simplex-scalar-replacement} with $g=f$ and permutation
invariance, $|\E\chi(U)-\mu|\le CN\varepsilon$. Together with the
variance bound, this gives, for $z=(\chi(U)-\mu)/b$,
\[
 \E z^2
 =\frac{\Var(\chi(U))+(\E\chi(U)-\mu)^2}{b^2}
 \le C/L_0^2.
\]

For $u\in H$, the directional derivative is $D_u=u^T\nabla_H$.
Gaussian convolution gives
\[
 \E f^{(3)}(\xi_\sigma)
 =\frac{\varphi_{1+\sigma^2}(1)}{1+\sigma^2}\ge c_0>0,
\]
where $c_0$ is universal since $\sigma^2\in[0,1]$. As
$|\varphi_1(z)-\varphi_1(0)|\le Cz^2$ and $f^{(3)}$ is bounded,
\eqref{eq:simplex-scalar-replacement} yields
\[
 \E[\varphi_1(z)f^{(3)}(U_1)]
 \ge \varphi_1(0)c_0-C(\varepsilon+L_0^{-2}).
\]
All mixed directional derivatives of $\chi$ through order four in
directions from $\{u_1,\ldots,u_N\}$ are bounded by universal
constants, since $\|u_j\|_1\le2$ and the derivatives of $f$ are
bounded. The chain rule and $b\ge1$ give
\[
 D_{u_j}^3h(U)
 =\frac{\varphi_1(z)}{b}\sum\nolimits_{k=1}^{N}(u_j)_k^3 f^{(3)}(U_k)
   +O(b^{-2}).
\]
By permutation invariance, $\E[\varphi_1(z)f^{(3)}(U_k)]$ is independent of $k$.
Also
$\sum_{k=1}^{N}(u_j)_k^3=(N-1)(N-2)/N^2\ge2/9$. For a sufficiently
large universal $L_0$ and then a sufficiently small universal
$\varepsilon_0$, these estimates yield, whenever $\varepsilon<\varepsilon_0$,
\[
 \E D_{u_j}^3h(U)\ge c_1/b,\qquad j=1,\ldots,N,
\]
with $c_1>0$ universal. The same chain rule bounds every mixed fourth
directional derivative of $h$ in directions from
$\{u_1,\ldots,u_N\}$ by $C/b$. In particular,
\[
 \sup_{x\in H}|D_{u_j}^4 h(x)|\le C/b,\qquad
 \sup_{x\in H}|\Delta_H^2h(x)|\le CN^2/b,
\]
because $\Delta_H=\sum_{j=1}^{N}D_{u_j}^2$.

In each replacement, the original term is treated by expanding
$h(U+X/\sqrt n)$ through order three. Conditional on $U$, each
displacement $a\sqrt{N/n}\,u_j$ has probability $a^{-2}/N$.
The absolute expected remainder is thus at most $Ca^2N^2/(bn^2)$.
For the Gaussian term, the function
$t\mapsto\E h(x+\sqrt t Z_P)$ has first and second derivatives
$\tfrac12\E\Delta_H h(x+\sqrt t Z_P)$ and
$\tfrac14\E\Delta_H^2 h(x+\sqrt t Z_P)$ for $t>0$, by Gaussian
integration by parts on $H$. These extend continuously to $t=0$
by boundedness and continuity of the derivatives. Taylor's formula
about $0$ at $t=1/n$ gives $h(x)+\Delta_Hh(x)/(2n)$ up to an error at most
$CN^2/(bn^2)$.

Conditional on $U$, both expansions contain
$h(U)+\Delta_Hh(U)/(2n)$, since both terms have covariance $P/n$.
The linear term for $X$ vanishes by centering.
The cubic contribution from the original term is at least
\[
 \frac{a\sqrt N}{6n^{3/2}}\sum\nolimits_{j=1}^{N}\E D_{u_j}^3h(U)
 \ge \frac{c_1aN^{3/2}}{6bn^{3/2}}.
\]
Summing the comparisons gives
\[
 \E h(S_n)-\E h(Z_P)
 \ge \frac{c_1aN^{3/2}}{6b\sqrt n}-\frac{Ca^2N^2}{bn}.
\]
The ratio of the error term to the first term is at most $C\varepsilon$.
Decreasing $\varepsilon_0$ if necessary proves the lower bound since
\[
 \frac{aN^{3/2}}{b\sqrt n}
 =\frac1{L_0}\frac{aN/\sqrt n}{1+aN/\sqrt n}
 \ge \frac1{2L_0}\min\{1,ad/\sqrt n\}.
\]
\end{proof}

\subsection{Proof of the lower bounds in Theorem~\ref{thm:moment-constraint}}
\label{sec:proof-moment-lower}
\begin{proof}
For $3\le q<4$, the relevant scale is
$L=(\widetilde{M}_q-d^{q/2})^{1/(q-2)}$. If $L\ge d$, the lattice law
has $M_q\le\widetilde{M}_q$ by \eqref{eq:lattice-moments}, and
\eqref{eq:lattice-lower} gives the required lower bound since $d+L\le2L$.
If $L<d$, the simplex law in Theorem~\ref{thm:simplex-lower}(i) with
$a=1$ suffices, since $M_q=d^{q/2}$ and $d+L<2d$.

For $q\ge4$, the simplex law in Theorem~\ref{thm:simplex-lower}(i)
with $a=(\widetilde{M}_q/d^{q/2})^{1/(q-2)}$ has
$M_q=\widetilde{M}_q$ and
\[
 ad=d^{(q-4)/(2(q-2))}\widetilde{M}_q^{1/(q-2)}.
\]
Theorem~\ref{thm:simplex-lower}(i) now gives the asserted lower bound.
\end{proof}

\section{Proof of Theorem~\ref{thm:independent-coordinates}}\label{sec:coordinate-proof}
\begin{proof}
Independence across the vectors and within each vector makes the $nd$
centered vectors $n^{-1/2}X_{ij}e_j$, $1\le i\le n$, $1\le j\le d$,
mutually independent. They sum to $S_n$ and their total covariance is
$\Var(S_n)=I_d$. The sum of their fourth moments is
$n^{-2}\sum_{i=1}^n\sum_{j=1}^d\E|X_{ij}|^4$, so Theorem~\ref{t}
gives the upper bound.

For the lower bound, all entries are independent copies of
$Y=3\sqrt2(2/3-T)$, where $T$ has density $2t\Ind_{(0,1)}(t)$.
Direct integration gives $\E Y=0$, $\E Y^2=1$, $\E Y^3=2\sqrt2/5$,
and $\E Y^4=12/5$. The vector $S_n^{[r]}$ is the projection of $S_n$
onto the first $r=\min\{d,n\}$ coordinates. This projection can only
decrease convex-set distance, since the inverse image of a convex
set is convex.

For a standard normal $\xi$, the function $f(x)=\E(x+\xi-1)_+$,
$x\in\R$, is convex and $1$-Lipschitz, with all positive-order
derivatives bounded. For $x\in\R^r$, the test function is given by
\[
 \chi(x)={\textstyle\sum_{j=1}^r f(x_j)},\qquad
 b=L_0\sqrt r,\qquad h(x)=\Phi\bigl((\chi(x)-r\E f(\xi))/b\bigr),
\]
where $L_0\ge1$ is a universal constant chosen below. The sublevel
sets of $h$ are convex. Since $0<h<1$, integration of their
probabilities gives
\[
 \begin{split}
 |\E h(S_n^{[r]})-\E h(Z_{I_r})|
 &\le \int_0^1
   |\P(h(S_n^{[r]})\le t)-\P(h(Z_{I_r})\le t)|\,dt\\
 &\le \dc(S_n^{[r]},Z_{I_r}).
 \end{split}
\]

Starting from the Gaussian sum, the $nr$ terms are replaced one at a
time by independent copies of $Y/\sqrt n$ in the corresponding
coordinates. At each replacement, $U$ is the sum with that term
omitted, independent of both candidate terms.
The coordinates of $U$ are independent and have variance at most one.
For $g=f$ or $g=f^{(3)}$, Taylor expansion through degree two
gives $|\E g(U_j)-\E g(\xi)|\le C/\sqrt n$ for each $j$.
Indeed, the first two moments agree in each scalar replacement and the
third-order remainder is at most $Cn^{-3/2}$ per replacement.
In the affected coordinate, adding the omitted centered Gaussian
term costs at most $\|g''\|_\infty/(2n)$. The other coordinates
already contain $n$ terms.
Independence and the Lipschitz bound give
$\Var(\chi(U))\le\sum_{j=1}^r\Var(U_j)\le r$.
Also $|\E\chi(U)-r\E f(\xi)|\le Cr/\sqrt n$. Thus, for
$z=(\chi(U)-r\E f(\xi))/b$,
\[
 \E z^2\le \frac{r+C^2r^2/n}{b^2}\le C/L_0^2 .
\]

Since $f''(x)=\varphi_1(x-1)$, Gaussian convolution gives
$\E f^{(3)}(\xi)=\varphi_2(1)/2>0$.
The chain rule gives
\[
 \partial_j^3 h(U)=b^{-1}\varphi_1(z) f^{(3)}(U_j)+R_j,\qquad
 |R_j|\le C/b^2,\qquad \sup_{x\in\R^r}|\partial_j^4 h(x)|\le C/b.
\]
The bound $|\varphi_1(z)-\varphi_1(0)|\le Cz^2$ and the scalar
replacement estimate give
\[
 b\,\E\partial_j^3 h(U)
 \ge \varphi_1(0)\varphi_2(1)/2
     -C\bigl(n^{-1/2}+L_0^{-2}+b^{-1}\bigr).
\]
For a sufficiently large universal $L_0$ and then a sufficiently large
universal $n_0$, the right-hand side is bounded below by a positive
universal constant for $n\ge n_0$.
Taylor expansion through degree three at $U$ now shows that
replacing a Gaussian term by the corresponding $Y/\sqrt n$ increases
the expectation by at least $c/(bn^{3/2})-C/(bn^2)$.
The first two moments match, the cubic term uses $\E Y^3>0$, and the
remainder bound follows from the fourth moments. Increasing $n_0$ if
necessary and summing over all $nr$ terms yields
\[
 \E h(S_n^{[r]})-\E h(Z_{I_r})\ge c\sqrt{r/n}
 =c\min\{1,\sqrt{d/n}\}.
\]
For $n<n_0$, the first coordinate of $S_n$ is bounded above by
$2\sqrt{2n_0}$. The halfspace
$A=\{x\in\R^d:x_1\le2\sqrt{2n_0}\}$ satisfies
$\P(S_n\in A)=1$ and
$\P(Z_{I_d}\in A)=\Phi(2\sqrt{2n_0})<1$.
Decreasing $c$ therefore covers these sample sizes.
\end{proof}

\section{Proof of the Gaussian perimeter bound}\label{sec:gaussian-perimeter-proof}
\begin{proof}[Proof of Lemma~\ref{lem:shell}]
For a bounded full-dimensional polytope $P$, the facets are
$F_1,\ldots,F_m$ and the corresponding unit outer normals are
$n_1,\ldots,n_m$. In this proof, $Q=\Sigma^{-1}$,
$h=\|Q\|_{\mathrm F}^{1/2}$ and $\psi=\varphi_\Sigma$.
The Gaussian mass of $P$ is $p=\P(Z_\Sigma\in P)$ and its Gaussian
perimeter is $S=\int_{\partial P}\psi\,dS$, where $dS$ is
$(d-1)$-dimensional Hausdorff measure. For the relative interior
$F_k^\circ$ of $F_k$, Nazarov's polar coordinate system is
$X_2(x,u)=x+un_k$, $x\in F_k^\circ$, $u>0$ \cite{nazarov}. The inequality $n_k^T(y-x)\le0$ gives
\[
 \|x+un_k-y\|^2\ge u^2+\|x-y\|^2,\qquad y\in P.
\]
Thus $x$ is the unique projection of $X_2(x,u)$ onto $P$.
The images of $F_k^\circ\times(0,\infty)$ under $X_2$ are disjoint
and lie outside $P$. On each image the volume element is $du\,dS(x)$.
Since $n_k^TQn_k\le h^2$,
\[
 \frac{\psi(x+un_k)}{\psi(x)}
 \ge e^{-u^2h^2/2}e^{-un_k^TQx}
 \ge e^{-u^2h^2/2}(1-un_k^TQx),\qquad u>0.
\]
Integrating over these images gives
\[
 1-p\ge\sqrt{\pi/2}\,\frac Sh-\frac{J}{h^2},\qquad
 J=\int_{\partial P}n(x)^TQx\,\psi(x)\,dS(x),
\]
where $n(x)=n_k$ on $F_k^\circ$. Facet boundaries are $dS$-null, and the
signed integrals are absolutely convergent since $P$ is bounded.
The divergence theorem and $\nabla\psi(x)=-Qx\psi(x)$ give
\[
 J=\E\!\left[\Ind_P(Z_\Sigma)
               \bigl(\tr Q-\|QZ_\Sigma\|^2\bigr)\right].
\]
Since $QZ_\Sigma\sim\mathcal N(0,Q)$, its squared norm has the law of $\sum_{i=1}^{d}\lambda_i N_i^2$, where
$\lambda_1,\ldots,\lambda_d$ are the eigenvalues of $Q$ and the $N_i$ are
i.i.d.\ $\mathcal N(0,1)$. Its mean is $\tr Q$ and its variance is
$2\|Q\|_{\mathrm F}^2=2h^4$. Centering the indicator and using
Cauchy--Schwarz gives $|J|\le h^2\sqrt{2p(1-p)}$. Hence
\[
 \frac Sh\le \sqrt{\frac2\pi}\bigl(1-p+\sqrt{2p(1-p)}\bigr)
 \le \frac{1+\sqrt3}{\sqrt{2\pi}}=C_{\scriptscriptstyle\mathrm{GSA}}.
\]
The last inequality follows from
$z+\sqrt2\sqrt{1-z^2}\le\sqrt3$ with $z=1-2p$.

For nonempty closed $C$, $C_t=C+tB_2^d$ because the distance to $C$
is attained. To bound the mass of $P+tB_2^d$, the supporting
hyperplanes of $P=\bigcap_{k=1}^{m}\{x:n_k^Tx\le a_k\}$ are moved
outwards by $u$, giving
\[
 f(x)=\max_{1\le k\le m}(n_k^Tx-a_k),\qquad x\in\R^d,
 \qquad D_u=\{x\in\R^d:f(x)\le u\},\qquad u\ge0.
\]
Each $D_u$ contains $P$ and is bounded. Indeed,
$\min_{\|v\|=1}\max_{1\le k\le m}n_k^Tv>0$, since a direction with
$n_k^Tv\le0$ for every $k$ would give a ray in $P$.
We have $P+tB_2^d\subseteq D_t$ because the normals have unit length.
The function $f$ is Lipschitz and $\|\nabla f\|=1$ almost everywhere,
since ties occur only on finitely many hyperplanes. Also
$\{f=u\}=\partial D_u$. Points with $f<u$ are interior. At a point
with $f=u$, moving in the direction $n_k$ for any active constraint
leaves $D_u$.
The coarea formula and the perimeter bound for polytopes give
\[
 \P(Z_\Sigma\in P+tB_2^d)-\P(Z_\Sigma\in P)
 \le\int_0^t\int_{\partial D_u}\psi(x)\,dS(x)\,du
 \le C_{\scriptscriptstyle\mathrm{GSA}}th,\qquad t\ge0.
\]

The extension to a nonempty compact convex set $C$ uses its support function
$h_C(v)=\sup_{x\in C}v^Tx$, $v\in\R^d$. A dense sequence $(u_j)_{j\ge1}$
on the unit sphere, starting with the $2d$ coordinate directions,
gives the polytopes
\[
 P_m=\bigcap\nolimits_{j=1}^{m}\{x\in\R^d:u_j^Tx\le h_C(u_j)+1/m\},
 \qquad m\ge2d.
\]
These nested polytopes lie in a common bounded box and contain
$C+m^{-1}B_2^d$. If $x\notin C$, strict separation gives a unit vector
$u$ with $u^Tx>h_C(u)$. By continuity, $u_j^Tx>h_C(u_j)$ for some $j$.
For all sufficiently large $m\ge j$, this strict inequality remains
true with $h_C(u_j)+1/m$ on the right, so $x\notin P_m$.
Thus $P_m\downarrow C$. If $x$ lies in every enlarged polytope, there
are points $y_m\in P_m$ with $\|x-y_m\|\le t$. A subsequential
limit belongs to $C$ and is within distance $t$ of $x$, proving
$P_m+tB_2^d\downarrow C+tB_2^d$. Continuity from above gives
\eqref{shell} for compact $C$.

For a nonempty unbounded closed convex $C$, the compact-set bound
applies to $C\cap mB_2^d$ for all sufficiently large $m$. These
sets and their $t$-enlargements increase to $C$ and $C_t$,
respectively, so continuity from below gives \eqref{shell}.

For any nonempty convex $C$, $C_t=(\overline C)_t$ and
$\gamma(C)=\gamma(\overline C)$. The latter equality follows from
the Lebesgue-null boundary in the full-dimensional case and from
containment in a proper affine subspace otherwise. Hence the bound
for $\overline C$ gives \eqref{shell}. The empty set has empty enlargement.
\end{proof}

\section{Proof of Lemma~\ref{lem:scalar-tilt}}\label{sec:scalar-tilt-proof}
\begin{proof}
Throughout the proof, $0<\beta\le r\le r_*$ and
$\lambda=\lambda_r(\beta)$. The integral $I(s)$ is defined in
Lemma~\ref{lem:scalar-tilt}. Its integrand, apart from the normalization
$1/\sqrt{2\pi}$, is continuous, equals one at $t=0$, and is bounded by
$e^{-\lambda t^2/2}$. It is positive near zero, and $\lambda>0$, so
$0<I(s)<\infty$ for $s\in\R$.

The functions are $h(s)=\log\cosh s$ and $g(s)=s^2/2-h(s)$,
$s\in\R$. For $s>0$, the series for $\operatorname{arctanh}u$ with
$u=\tanh s$ gives
\[
 \frac{g'(s)}{h'(s)}=\frac{s-u}{u}
 =\sum\nolimits_{k=1}^{\infty}\frac{u^{2k}}{2k+1}
 \ge-\frac13\log(1-u^2)=\frac23 h(s).
\]
Since $g(0)=h(0)=0$, integration of $g'\ge(2/3)hh'$ and evenness give
$g(s)\ge h(s)^2/3$. Thus
\begin{equation}\label{eq:app-a-gap}
 b h(s)-g(s)\le 3b^2/4,\qquad b\ge0,\quad s\in\R.
\end{equation}

For $t_k=k\pi/\sqrt\beta$, the contribution of the $k$th interval is
\[
 I_k(s)=\frac1{\sqrt{2\pi}}
 \int_{|t-t_k|\le\pi/(2\sqrt\beta)} e^{-\lambda t^2/2}
 \bigl(1-\operatorname{sech}^2s\,\sin^2(\sqrt\beta\,t)\bigr)^{1/(2\beta)}\,dt,
 \qquad s\in\R,\quad k\in\mathbb Z.
\]
These intervals cover $\R$ with disjoint interiors, so
$I(s)=\sum_{k\in\mathbb Z}I_k(s)$.
With $a=4/(\pi^2\cosh^2s)$ and $z=t-t_k$, the inequalities
$|\sin v|\ge2|v|/\pi$ for $|v|\le\pi/2$ and $1-x\le e^{-x}$
give, on this interval,
\begin{equation}\label{eq:app-a-chord}
 \bigl(1-\operatorname{sech}^2 s\,\sin^2(\sqrt\beta\,t)\bigr)^{1/(2\beta)}
 \le e^{-az^2/2}.
\end{equation}
Also $\lambda/\beta\ge A_0$, and $\log x\le x-1$ with $L\le A_0$ gives
$\lambda\le A_0\max\{\beta,r^2\}\le D:=A_0r_*$.

When $\cosh^2 s\le2$, the estimate uses
$T^2=(\lambda+\operatorname{sech}^2 s)/\beta$. The bound
$\sin^2(\sqrt\beta\,t)\le\beta t^2$ in the Gaussian factor and the
substitution $y=T\sin(\sqrt\beta\,t)$ on the positive half of the
central interval give
\[
 \sqrt{\lambda+\operatorname{sech}^2s}\, I_0(s)
 \le \sqrt{2/\pi}\int_0^T e^{-y^2/2}(1-y^2/T^2)^{-1/2}\,dy
 \le 1+4T^{-2}.
\]
The estimate on $[0,T/\sqrt2]$ follows from
$(1-w)^{-1/2}\le1+w$ for $0\le w\le1/2$ and the Gaussian second moment.
On $[T/\sqrt2,T]$ the contribution is at most
\[
 \sqrt{2/\pi}(\pi T/4)e^{-T^2/4}<3T^{-2},
\]
since $\sup_{T>0}T^3 e^{-T^2/4}=6^{3/2}e^{-3/2}$.
As $T^{-2}\le2\beta$ and
$(1+\lambda)/(\lambda+\operatorname{sech}^2 s)\le\cosh^2 s$, we obtain
\begin{equation}\label{eq:app-a-central}
 I_0(s)\sqrt{1+\lambda}\le\cosh s\,(1+8\beta).
\end{equation}

Still assuming $\cosh^2s\le2$, the function is
$c(x)=P_0/(1+P_0x)$, $x\ge0$, where $P_0=\pi^2/2$.
Combining \eqref{eq:app-a-chord} with the Gaussian factor, completing
the square, and extending the integral to $\R$ give
$I_k(s)\le(a+\lambda)^{-1/2}\exp[-a\lambda t_k^2/(2(a+\lambda))]$.
Since $1/P_0\le a\le1$, normalization by $\sqrt{1+\lambda}/\cosh s$
bounds the prefactor by $\pi/2$ and the exponential by
$e^{-c(\lambda)(\lambda/\beta)k^2}$. Summing with
$k^2\ge1+3(k-1)$ for $k\ge1$ gives
\begin{equation}\label{eq:app-a-outer}
 \frac{\sqrt{1+\lambda}}{\cosh s}
 \sum\nolimits_{\substack{k\in\mathbb Z\\ k\ne0}} I_k(s)
 \le \frac{\pi\,e^{-c(\lambda)\lambda/\beta}}{1-e^{-3A_0 c(D)}}
 <3.2\,e^{-c(\lambda)\lambda/\beta}.
\end{equation}
For the stated parameters, $3A_0c(D)>7$ and
$\pi/(1-e^{-7})<3.2$.

If $\beta\le r^2$, then $\lambda\le A_0r_*^2$, so
$3.2e^{-A_0c(\lambda)}<1/16$ and $Lc(\lambda)>1$. Consequently
\[
 3.2e^{-c(\lambda)\lambda/\beta}
 =3.2e^{-A_0c(\lambda)}(\beta/r^2)^{Lc(\lambda)}
 <{\beta}/{(16r^2)}.
\]
If $\beta\ge r^2$, the function $F(x)=e^{-\omega/(1+\omega x)}$,
$0\le x\le r_*$, with $\omega=A_0P_0$, is convex because $F''(x)$ has the sign of
$\omega-2(1+\omega x)>0$ on this interval, since
$\omega(1-2r_*)>2$.
The preceding case at $\beta=r^2$ gives $3.2F(r^2)<1/16$, and
$3.2r F(r)\le 3.2r_*e^{-A_0c(D)}<1/16$.
By convexity and the two endpoint bounds,
$3.2F(x)\le x/(16r^2)$ for $r^2\le x\le r$.
In particular, $3.2F(\beta)\le\beta/(16r^2)$, and in this case
$e^{-c(\lambda)\lambda/\beta}=F(\beta)$.

Combining \eqref{eq:app-a-central} and \eqref{eq:app-a-outer} and using
$\log(1+x)\le x$ and \eqref{eq:app-a-gap}, we obtain
\[
 \beta\log\!\bigl(I(s)\sqrt{1+\lambda}\bigr)-g(s)
 \le \beta h(s)-g(s)+8\beta^2+\frac{\beta^2}{16r^2}
 \le \frac{35}{4}\beta^2+\frac{\beta^2}{16r^2}.
\]

When $\cosh^2 s>2$, the Gaussian factor in the central interval
is at most $1$.
On each outer interval, the factor $e^{-\lambda t^2/2}$ is at most
$e^{-\pi^2A_0(|k|-1/2)^2/2}$, since
$|t|\ge (|k|-1/2)\pi/\sqrt\beta$ and $\lambda/\beta\ge A_0$.
Integrating \eqref{eq:app-a-chord} over the line and summing with
$(k-1/2)^2\ge1/4+2(k-1)$ for $k\ge1$ gives
\[
 \frac{I(s)\sqrt{1+\lambda}}{\cosh s}
 \le\frac{\pi\sqrt{1+D}}{2}
 \left(1+\frac{2e^{-\pi^2 A_0/8}}{1-e^{-\pi^2 A_0}}\right)
 <e.
\]
Since $h(s)>\tfrac12\log2>1/3$, we have $1+h(s)\le4h(s)$. Therefore
\[
 \beta\log\!\bigl(I(s)\sqrt{1+\lambda}\bigr)-g(s)
 \le 4\beta h(s)-g(s)\le12\beta^2
\]
by \eqref{eq:app-a-gap}. The two cases give
\eqref{eq:app-a-scalar} by the definition of $K(r)$.
\end{proof}

\paragraph*{Acknowledgement of AI use:} The author used AI tools for improving results in the initial versions of the paper, LaTeX formatting, assistance with proof exposition and organization, and for language editing. The author reviewed and revised the manuscript, verified the mathematical arguments, and takes full responsibility for its content.


\end{document}